# A new representation of the Stieltjes constants

Donal F. Connon

dconnon@btopenworld.com

11 January 2022

**Abstract**

Inter alia, this paper shows how a new representation of the Stieltjes constants $\gamma_n(x)$ leads to the evaluation of $\varsigma^{(n+1)}(0,x)$. For example, we show that for $0 < x < 1$

$$\gamma_1(x) = \lim_{s \to 0} \gamma_1(x,s)$$

where $\gamma_1(x,s)$ is defined by

$$\gamma_1(x,s) = \sum_{n=1}^{\infty} \frac{[\gamma + \log(2\pi n)]^2 \cos(2\pi nx)}{n^s} - \frac{1}{2}\varsigma(2)\sum_{n=1}^{\infty} \frac{\cos(2\pi nx)}{n^s}$$

$$-\pi \sum_{n=1}^{\infty} \frac{[\gamma + \log(2\pi n)]\sin(2\pi nx)}{n^s}$$

Using this, we also show that for $0 < x < 1$

$$\varsigma''(0,x) = \sum_{n=1}^{\infty} \frac{[\gamma + \log(2\pi n)]^2 \sin(2\pi nx)}{n\pi} - \frac{1}{4}\varsigma(2)\sum_{n=1}^{\infty} \frac{\sin(2\pi nx)}{n\pi}$$

$$+\sum_{n=1}^{\infty} \frac{[\gamma + \log(2\pi n)]\cos(2\pi nx)}{n}$$

## 1. Introduction

This paper may usefully be read in conjunction with the author's earlier related paper, *A Ramanujan enigma involving the first Stieltjes constant* [27], which was originally published in 2019 (and recently amended). Inter alia, this paper showed that for $0 < x < 1$

(1.1)  $\gamma_1(1-x) - \gamma_1(x) = \pi[\gamma + \log(2\pi)]\cot(\pi x) + 2\pi \lim_{s \to 1} \sum_{n=1}^{\infty} \frac{\sin(2n\pi x)\log n}{n^{1-s}}$

where $\gamma_n(x)$, the generalised Stieltjes constants, are the coefficients in the Laurent expansion of the Hurwitz zeta function $\varsigma(s,x)$ about $s = 1$

$$(1.2.1) \qquad \varsigma(s,x) = \sum_{n=0}^{\infty} \frac{1}{(n+x)^s}$$

$$(1.2.2) \qquad \varsigma(s,x) = \frac{1}{s-1} + \sum_{n=0}^{\infty} \frac{(-1)^n}{n!} \gamma_n(x)(s-1)^n$$

We have

$$(1.3) \qquad \gamma_0(x) = -\psi(x)$$

where $\psi(x)$ is the digamma function which is the logarithmic derivative of the gamma function $\psi(x) = \frac{d}{du} \log \Gamma(x)$. It is easily seen from the definition of the Hurwitz zeta function that $\varsigma(s,1) = \varsigma(s)$ and accordingly that $\gamma_n(1) = \gamma_n$.

The analysis presented in this paper is somewhat Eulerian in nature but, notwithstanding the apparent lack of rigour, it does generate valid results. Unless otherwise specified, most of the formulae in this paper apply in the open interval $(0,1)$.

We showed in [27] that for $0 < x < 1$

$$(1.4) \qquad \lim_{s \to 0} \sum_{n=1}^{\infty} \frac{\cos(2n\pi x)}{n^s} = -\frac{1}{2}$$

$$(1.5) \qquad \lim_{s \to 0} \sum_{n=1}^{\infty} \frac{\sin(2n\pi x)}{n^s} = \frac{1}{2} \cot \pi x$$

Srivastava and Tsumura [45, p.293] reported for $\text{Re}(s) > 1$

$$(1.6) \qquad \sum_{n=1}^{\infty} \frac{\cos(n\pi/3)}{n^s} = \frac{1}{2}(6^{1-s} - 3^{1-s} - 2^{1-s} + 1)\varsigma(s)$$

$$(1.7) \qquad \sum_{n=1}^{\infty} \frac{\cos(2n\pi/3)}{n^s} = \frac{1}{2}(3^{1-s} - 1)\varsigma(s)$$

$$(1.8) \qquad \sum_{n=1}^{\infty} \frac{\cos(n\pi/2)}{n^s} = 2^{-s}(2^{1-s} - 1)\varsigma(s)$$

and it is unlikely to be just an interesting coincidence that these identities also concur with (1.4) in the limit $s \to 0$. It was also reported by Srivastava and Tsumura [45, p.293] that for $\text{Re}(s) > 1$

$$(1.9) \qquad \sum_{n=1}^{\infty} \frac{\sin(n\pi/3)}{n^s} = \sqrt{3} \left\{ \frac{3^{-s} - 1}{2} \varsigma(s) + 6^{-s} \left[ \varsigma\left(s, \frac{1}{6}\right) + \varsigma\left(s, \frac{1}{3}\right) \right] \right\}$$

(1.10) $$\sum_{n=1}^{\infty} \frac{\sin(2n\pi/3)}{n^s} = \sqrt{3}\left\{\frac{3^{-s}-1}{2}\varsigma(s) + 3^{-s}\varsigma\left(s,\frac{1}{3}\right)\right\}$$

(1.11) $$\sum_{n=1}^{\infty} \frac{\sin(n\pi/2)}{n^s} = (2^{-s}-1)\varsigma(s) + 2^{1-2s}\varsigma\left(s,\frac{1}{4}\right)$$

and it may be noted that these identities also concur with (1.5) in the limit $s \to 0$.

The Lerch zeta function $L(x,a,s)$ is defined [38] for $\mathrm{Re}(s) > 1$ by

$$L(x,a,s) = \sum_{n=0}^{\infty} \frac{e^{2\pi i n x}}{(a+n)^s}$$

and this may be analytically continued to the whole $s$-plane by means of a contour integral.

We see that

$$\sum_{n=0}^{\infty} \frac{e^{2\pi i n x}}{(a+n)^s} = \frac{1}{a^s} + \sum_{n=1}^{\infty} \frac{e^{2\pi i n x}}{(a+n)^s}$$

Apostol [1] proved that

$$L(x,a,0) = \frac{1}{2} + \frac{1}{2}i\cot\pi x$$

which *suggests* that

$$\lim_{s \to 0} \sum_{n=1}^{\infty} \frac{e^{2\pi i n x}}{(a+n)^s} = -\frac{1}{2} + \frac{1}{2}i\cot\pi x$$

In a slightly different approach, we see that $L(x,a,0)$ is independent of the value of $a$, and with $a = 1$ we obtain

$$L(x,1,0) = \frac{1}{2} + \frac{1}{2}i\cot\pi x$$

We also see that

$$L(x,1,s) = \sum_{n=0}^{\infty} \frac{e^{2\pi i n x}}{(1+n)^s}$$

$$= \sum_{n=1}^{\infty} \frac{e^{2\pi i (n-1) x}}{n^s}$$

$$= e^{-2\pi ix} \sum_{n=1}^{\infty} \frac{e^{2\pi inx}}{n^s}$$

as noted in [38, p.39]. Therefore we have

$$\lim_{s \to 0} \sum_{n=1}^{\infty} \frac{e^{2\pi inx}}{n^s} = e^{2\pi ix}\left[\frac{1}{2} + \frac{1}{2}i\cot\pi x\right]$$

$$= -\frac{1}{2} + \frac{1}{2}i\cot\pi x$$

The following identities are recorded by Hansen [50, pp. 223 & 244] for $\operatorname{Re}(s) > 1$ and $0 < x < 2\pi$

(1.12)

$$\sum_{n=1}^{\infty} \frac{\sin(nx+y)}{n^s} = \frac{(2\pi)^s}{2\Gamma(s)}\operatorname{cosec}(\pi s)\left[\cos\left(y - \frac{\pi s}{2}\right)\varsigma\left(1-s, \frac{x}{2\pi}\right) - \cos\left(y + \frac{\pi s}{2}\right)\varsigma\left(1-s, 1-\frac{x}{2\pi}\right)\right]$$

(1.13)

$$\sum_{n=1}^{\infty} \frac{\cos(nx+y)}{n^s} = \frac{(2\pi)^s}{2\Gamma(s)}\operatorname{cosec}(\pi s)\left[\sin\left(y + \frac{\pi s}{2}\right)\varsigma\left(1-s, 1-\frac{x}{2\pi}\right) - \sin\left(y - \frac{\pi s}{2}\right)\varsigma\left(1-s, \frac{x}{2\pi}\right)\right]$$

and we may note that the second identity may be obtained by substituting $y \to y + \frac{\pi}{2}$ (or by differentiating the first one with respect to $y$).

We write (1.12) as

$$\sum_{n=1}^{\infty} \frac{\sin(2\pi nx + y)}{(2\pi n)^s} = \frac{1}{2\Gamma(s)}\operatorname{cosec}(\pi s)\left[\cos\left(y - \frac{\pi s}{2}\right)\varsigma(1-s, x) - \cos\left(y + \frac{\pi s}{2}\right)\varsigma(1-s, 1-x)\right]$$

and note that

$$\cos\left(y - \frac{\pi s}{2}\right)\varsigma(1-s, x) - \cos\left(y + \frac{\pi s}{2}\right)\varsigma(1-s, 1-x)$$

$$= \left[\cos\left(y - \frac{\pi s}{2}\right) - \cos y\right]\varsigma(1-s, x) - \left[\cos\left(y + \frac{\pi s}{2}\right) - \cos y\right]\varsigma(1-s, 1-x)$$

$$+ \left[\varsigma(1-s, x) - \varsigma(1-s, 1-x)\right]\cos y$$

$$= \frac{\cos\left(y - \frac{\pi s}{2}\right) - \cos y}{s} s\varsigma(1-s,x) - \frac{\cos\left(y + \frac{\pi s}{2}\right) - \cos y}{s} s\varsigma(1-s,1-x)$$

$$+ \left[\varsigma(1-s,x) - \varsigma(1-s,1-x)\right]\cos y$$

We determine the limit as $s \to 0$ using L'Hôpital's rule and

$$\frac{1}{\Gamma(s)}\operatorname{cosec}(\pi s) = \frac{1}{\pi}\Gamma(1-s)$$

to obtain

$$\lim_{s \to 0} \sum_{n=1}^{\infty} \frac{\sin(2\pi nx + y)}{(2\pi n)^s} = -\frac{1}{2}\sin y + \frac{1}{2\pi}\cos y\left[\gamma_0(x) - \gamma_0(1-x)\right]$$

This may be written as

$$\lim_{s \to 0} \sum_{n=1}^{\infty} \frac{\sin(2\pi nx + y)}{(2\pi n)^s} = -\frac{1}{2}\sin y - \frac{1}{2\pi}\cos y\left[\psi(x) - \psi(1-x)\right]$$

or equivalently for $0 < x < 1$

(1.14) $$\lim_{s \to 0} \sum_{n=1}^{\infty} \frac{\sin(2\pi nx + y)}{(2\pi n)^s} = -\frac{1}{2}\sin y + \frac{1}{2}\cos y \cot(\pi x)$$

$$= \frac{1}{2}\frac{\cos(\pi x + y)}{\sin(\pi x)}$$

With $y = \frac{\pi}{2}$ we obtain (1.4) and $y = 0$ gives us (1.5).

Formally differentiating (1.14) gives us for $0 < x < 1$

(1.15) $$\lim_{s \to 0} \sum_{n=1}^{\infty} \frac{\cos(2\pi nx + y)}{(2\pi n)^s} = -\frac{1}{2}\frac{\sin(\pi x + y)}{\sin(\pi x)}$$

With $y = 0$ and $x = \frac{1}{2}$ in (1.15) we obtain

(1.16) $$\lim_{s \to 0} \sum_{n=1}^{\infty} \frac{(-1)^n}{(2\pi n)^s} = -\frac{1}{2}$$

Noting that $\sum_{n=1}^{\infty} \frac{(-1)^n}{n^s} = -(1 - 2^{1-s})\varsigma(s)$ we see that $\lim_{s \to 0} \sum_{n=1}^{\infty} \frac{(-1)^n}{n^s} = \varsigma(0) = -\frac{1}{2}$ (see [51]).

Reference may also be profitably made to [16] and [35, p.101].

It therefore appears that by using analytic continuation, or otherwise, it may be possible to justify the apparently unorthodox mathematics contained in this paper. In this regard, I quote from Hardy's book, *Divergent Series* [34, p.5]:

"It is natural to suppose that the other formulae will prove to be correct, and our transformations justifiable, if they are interpreted appropriately. We should then be able to regard the transformations as shorthand representations of more complicated processes justifiable by the ordinary canons of analysis."

Further important elucidation may be obtained from Section 9 below which sets out some pertinent analysis which was kindly provided to me by Bernard Candelpergher.

## 2. An appearance of the digamma function

We have the well-known Hurwitz's formula for the Fourier series of the Riemann zeta function $\varsigma(s,x)$ as reported in Titchmarsh's treatise [46, p.37]

$$(2.1) \quad \varsigma(s,x) = 2\Gamma(1-s)\left[\sin\left(\frac{\pi s}{2}\right)\sum_{n=1}^{\infty}\frac{\cos 2\pi nx}{(2\pi n)^{1-s}} + \cos\left(\frac{\pi s}{2}\right)\sum_{n=1}^{\infty}\frac{\sin 2\pi nx}{(2\pi n)^{1-s}}\right]$$

where $\text{Re}(s) < 0$ and $0 < x \leq 1$. In 2000, Boudjelkha [12] showed that this formula also applies in the region $0 < \text{Re}(s) < 1$ with $0 < x < 1$. It may be noted that (2.1) reduces to Riemann's functional equation for $\varsigma(s)$ when $x = 1$.

It is curious that neither Apostol [4, p.257], Titchmarsh [46, p.37] nor Whittaker and Watson [49, p.268] specifically highlight the rather obvious fact that (2.1) remains valid in the case $s = 0$ when $0 < x < 1$. We then have

$$\varsigma(0,x) = \sum_{n=1}^{\infty}\frac{\cos 2\pi nx}{\pi n}$$

and it is well known that for $0 < x < 1$

$$\frac{1}{2} - x = \sum_{n=1}^{\infty}\frac{\cos 2\pi nx}{\pi n}$$

We have [4, p.264] for $m \in N_0$

$$\varsigma(-m,x) = -\frac{B_{m+1}(x)}{m+1}$$

where $B_m(x)$ are the Bernoulli polynomials. In particular, we have

$$B_1(x) = x - \frac{1}{2}$$

References to the first and second derivatives of (2.1) computed at $s = 0$ are contained in Sections 5 and 8 below. □

Since
$$\cos\left(\tfrac{\pi}{2}s\right)\varsigma(s,x) = \sin\left(\tfrac{\pi}{2}(1-s)\right)\varsigma(s,x)$$

we have
$$\lim_{s\to 1}\cos\left(\tfrac{\pi}{2}s\right)\varsigma(s,x) = \lim_{s\to 1}\frac{\sin\left(\tfrac{\pi}{2}(1-s)\right)}{\tfrac{\pi}{2}(1-s)}\tfrac{\pi}{2}(1-s)\varsigma(s,x)$$

$$= \lim_{s\to 1}\tfrac{\pi}{2}(1-s)\varsigma(s,x) = -\tfrac{\pi}{2}$$

and hence we have
$$\lim_{s\to 1}\cos\left(\tfrac{\pi}{2}s\right)\varsigma(s,x) = -\tfrac{\pi}{2}$$

We now multiply the Hurwitz Fourier series (2.1) by $\cos\left(\tfrac{\pi}{2}s\right)$ and consider the limit as $s \to 0$. We see using the right-hand side of (2.1) that (where we initially just consider the first summation)

$$2\Gamma(1-s)\cos\left(\tfrac{\pi}{2}s\right)\sin\left(\tfrac{\pi}{2}s\right)\sum_{n=1}^{\infty}\frac{\cos 2n\pi x}{(2\pi n)^{1-s}} = \Gamma(1-s)\sin(\pi s)\sum_{n=1}^{\infty}\frac{\cos 2n\pi x}{(2\pi n)^{1-s}}$$

$$= (1-s)\Gamma(1-s)\frac{\sin(\pi s)}{1-s}\sum_{n=1}^{\infty}\frac{\cos 2n\pi x}{(2\pi n)^{1-s}}$$

$$= \Gamma(2-s)\frac{\sin(\pi s)}{1-s}\sum_{n=1}^{\infty}\frac{\cos 2n\pi x}{(2\pi n)^{1-s}}$$

Therefore we have

$$\lim_{s\to 1}2\Gamma(1-s)\cos\left(\tfrac{\pi}{2}s\right)\sin\left(\tfrac{\pi}{2}s\right)\sum_{n=1}^{\infty}\frac{\cos 2n\pi x}{(2\pi n)^{1-s}} = -\pi\lim_{s\to 1}\sum_{n=1}^{\infty}\frac{\cos 2n\pi x}{(2\pi n)^{1-s}}$$

and employing (1.4) we see that

$$\lim_{s\to 1}2\Gamma(1-s)\cos\left(\tfrac{\pi}{2}s\right)\sin\left(\tfrac{\pi}{2}s\right)\sum_{n=1}^{\infty}\frac{\cos 2n\pi x}{(2\pi n)^{1-s}} = -\frac{\pi}{2}$$

We also have

$$\lim_{s\to 1}2\Gamma(1-s)\cos^2\left(\tfrac{\pi}{2}s\right)\sum_{n=1}^{\infty}\frac{\sin 2n\pi x}{(2\pi n)^{1-s}} = \lim_{s\to 1}2[(1-s)\Gamma(1-s)]\frac{\cos\left(\tfrac{\pi}{2}s\right)}{1-s}\cos\left(\tfrac{\pi}{2}s\right)\sum_{n=1}^{\infty}\frac{\sin 2n\pi x}{(2\pi n)^{1-s}}$$

and using (1.5) we see that this vanishes.
$\square$

We may express (2.1) more concisely as

$$\text{(2.2)} \qquad \varsigma(s,x) = 2\Gamma(1-s)\sum_{n=1}^{\infty} \frac{\sin\left(2n\pi x + \frac{\pi s}{2}\right)}{(2\pi n)^{1-s}}$$

and multiplying by $s-1$ results in

$$\text{(2.3)} \qquad (s-1)\varsigma(s,x) = -2\Gamma(2-s)\sum_{n=1}^{\infty} \frac{\sin\left(2n\pi x + \frac{\pi s}{2}\right)}{(2\pi n)^{1-s}}$$

The following limit is easily deduced

$$\lim_{s\to 1}(s-1)\varsigma(s,x) = -2\lim_{s\to 1}\sum_{n=1}^{\infty} \frac{\sin\left(2n\pi x + \frac{\pi s}{2}\right)}{(2\pi n)^{1-s}}$$

$$= -2\lim_{s\to 1}\cos\left(\tfrac{\pi s}{2}\right)\sum_{n=1}^{\infty} \frac{\sin(2n\pi x)}{(2\pi n)^{1-s}} - 2\lim_{s\to 1}\sin\left(\tfrac{\pi s}{2}\right)\sum_{n=1}^{\infty} \frac{\cos(2n\pi x)}{(2\pi n)^{1-s}}$$

and, since we showed in [27] that the limit $\displaystyle\lim_{s\to 1}\sum_{n=1}^{\infty} \frac{\sin(2n\pi x)}{(2\pi n)^{1-s}}$ is finite (see (2.7) below), we deduce that

$$\lim_{s\to 1}(s-1)\varsigma(s,x) = -2\lim_{s\to 1}\sum_{n=1}^{\infty} \frac{\cos(2n\pi x)}{(2\pi n)^{1-s}}$$

We see from (1.2) that $\displaystyle\lim_{s\to 1}(s-1)\varsigma(s,x) = 1$ and thus we have the important limit (see also [27] for further discussion)

$$\text{(2.4)} \qquad \lim_{s\to 1}\sum_{n=1}^{\infty} \frac{\cos(2n\pi x)}{(2\pi n)^{1-s}} = -\frac{1}{2}$$

Differentiating (2.3) gives us

$$\text{(2.5)} \qquad \frac{\partial}{\partial s}(s-1)[\varsigma(s,x)] = -2\Gamma(2-s)\sum_{n=1}^{\infty} \frac{\log(2\pi n)\sin\left(2n\pi x + \frac{\pi s}{2}\right) + \frac{\pi}{2}\cos\left(2n\pi x + \frac{\pi s}{2}\right)}{(2\pi n)^{1-s}}$$

$$+ 2\Gamma'(2-s)\sum_{n=1}^{\infty} \frac{\sin\left(2n\pi x + \frac{\pi s}{2}\right)}{(2\pi n)^{1-s}}$$

Since $\left.\dfrac{\partial}{\partial s}(s-1)\varsigma(s,x)\right|_{s=1} = \gamma_0(x)$ we see from (1.3) that

$$\psi(x) = \lim_{s\to 1}\sum_{n=1}^{\infty} \frac{2\log(2\pi n)\cos(2n\pi x) - \pi\sin(2n\pi x) + 2\gamma\cos(2n\pi x)}{(2\pi n)^{1-s}}$$

where we have used $\Gamma'(1) = -\gamma$. It should be noted that in determining the above limit we have specifically used (1.1) to see that $\lim_{s \to 1} \sum_{n=1}^{\infty} \frac{\sin(2n\pi x)\log n}{n^{1-s}}$ is finite and accordingly determine that

$$\lim_{s \to 1} \cos\left(\tfrac{\pi s}{2}\right) \sum_{n=1}^{\infty} \frac{\sin(2n\pi x)\log n}{n^{1-s}} = 0$$

We may express $\psi(x)$ as

(2.6) $$\psi(x) = \lim_{s \to 1} \sum_{n=1}^{\infty} \frac{2[\gamma + \log(2\pi n)]\cos(2n\pi x) - \pi \sin(2n\pi x)}{n^{1-s}}$$

We showed in [27] that

(2.7) $$\lim_{s \to 1} \sum_{n=1}^{\infty} \frac{\sin 2n\pi x}{n^{1-s}} = \frac{1}{2}\cot(\pi x)$$

and hence, using (2.4) and (2.7), we may express (2.6) as

(2.8) $$\psi(x) = -[\gamma + \log(2\pi)] - \frac{\pi}{2}\cot(\pi x) + \lim_{s \to 1} 2\sum_{n=1}^{\infty} \frac{\log n \cos(2n\pi x)}{n^{1-s}}$$

or equivalently

(2.8.1) $$\psi(x) = -\frac{\pi}{2}\cot(\pi x) + \lim_{s \to 1} 2\sum_{n=1}^{\infty} \frac{[\gamma + \log(2\pi n)]\cos(2n\pi x)}{n^{1-s}}$$

$\square$

With $x = \tfrac{1}{2}$ in (2.8.1) we have

$$\psi(\tfrac{1}{2}) = \lim_{s \to 1} 2\sum_{n=1}^{\infty} \frac{(-1)^n [\gamma + \log(2\pi n)]}{n^{1-s}}$$

$$= -\lim_{s \to 1} 2[\gamma + \log(2\pi)]\varsigma_a(1-s) \sum_{n=1}^{\infty} \frac{(-1)^n}{n^{1-s}} + \lim_{s \to 1} 2\varsigma_a'(1-s)$$

where $\varsigma_a(s) = \sum_{n=1}^{\infty} \frac{(-1)^{n+1}}{n^s}$. It is known that

$$\varsigma_a(s) = [1 - 2^{1-s}]\varsigma(s) \text{ for } \operatorname{Re}(s) > 0; \ s \neq 1$$

Therefore

$$\varsigma_a(0) = -\varsigma(0)$$

and hence

$$\varsigma_a(0) = \tfrac{1}{2}$$

Accordingly, this 'corresponds' with

$$\lim_{s \to 1} \sum_{n=1}^{\infty} \frac{(-1)^n}{n^{1-s}} = -\frac{1}{2}$$

which concurs with (1.4) when $x = \tfrac{1}{2}$.

We have the derivative

$$\varsigma_a'(s) = [1 - 2^{1-s}]\varsigma'(s) + 2^{1-s} \log 2 \, \varsigma(s)$$

and hence

$$\varsigma_a'(0) = -\varsigma'(0) - \log 2$$

$$= \frac{1}{2} \log(2\pi) - \log 2$$

Therefore we see that

$$\varsigma_a'(0) = \frac{1}{2} \log\left(\frac{\pi}{2}\right)$$

Accordingly, this 'corresponds' with

$$\lim_{s \to 1} \sum_{n=1}^{\infty} \frac{(-1)^n \log n}{n^{1-s}} = \frac{1}{2} \log\left(\frac{\pi}{2}\right)$$

which corresponds with equation (4.38) in Candelpergher's paper [14, p.139] for the 'Euler summation'

$$\sum_{n \geq 1}^{E} (-1)^n \log n = \frac{1}{2} \log\left(\frac{\pi}{2}\right)$$

Hence. we formally obtain

$$\psi(\tfrac{1}{2}) = -2[\gamma + \log(2\pi)]\varsigma_a(0) + 2\varsigma_a'(0)$$

$$= -[\gamma + \log(2\pi)] + \log(\pi/2)$$

and this concurs with the known value [44, p.20]

$$\psi(\tfrac{1}{2}) = -\gamma - 2\log 2$$

□

Letting $x \to 1-x$ in (2.6) gives us

$$\psi(1-x) = \lim_{s \to 1} \sum_{n=1}^{\infty} \frac{2[\gamma + \log(2\pi n)]\cos(2n\pi x) + \pi \sin(2n\pi x)}{n^{1-s}}$$

and hence we obtain

$$\psi(x) - \psi(1-x) = -2\pi \lim_{s \to 1} \sum_{n=1}^{\infty} \frac{\sin(2n\pi x)}{n^{1-s}}$$

Using (2.7) we then obtain the well-known functional equation

$$\psi(x) - \psi(1-x) = -\pi \cot(\pi x)$$

□

We multiply (2.8) by $\sin \pi x$ to obtain

(2.9) $$\psi(x)\sin \pi x = -[\gamma + \log(2\pi)]\sin \pi x - \frac{\pi}{2}\cos \pi x + \lim_{s \to 1} 2\sum_{n=1}^{\infty} \frac{\log n \cos(2n\pi x)\sin \pi x}{n^{1-s}}$$

and we immediately note the close structural similarity to Lerch's trigonometric series expansion [39] for the digamma function for $0 < x < 1$ (see for example, Gronwall's paper [32, p.105] and Nielsen's book [42, p.204]

(2.10) $$\psi(x)\sin \pi x = -[\gamma + \log(2\pi)]\sin \pi x - \frac{\pi}{2}\cos \pi x - \sum_{n=1}^{\infty} \sin(2n+1)\pi x \cdot \log\left(1 + \frac{1}{n}\right)$$

Equating (2.9) and (2.10) gives us

(2.11) $$-\sum_{n=1}^{\infty} \sin(2n+1)\pi x \cdot \log\left(1 + \frac{1}{n}\right) = \lim_{s \to 1} 2\sum_{n=1}^{\infty} \frac{\log n \cos(2n\pi x)\sin \pi x}{n^{1-s}}$$

We now integrate (2.11) over the interval $[\tfrac{1}{2}, u]$ with $0 < u < 1$ where we assume that

$$\int_{\frac{1}{2}}^{u} \lim_{s \to 1} 2\sum_{n=1}^{\infty} \frac{\log n \cos(2n\pi x)\sin \pi x}{n^{1-s}} dx = \lim_{s \to 1} \sum_{n=1}^{\infty} \int_{\frac{1}{2}}^{u} \frac{\log n \cos(2n\pi x)\sin \pi x}{n^{1-s}} dx$$

Using

$$\int \cos(2n\pi x)\sin \pi x \, dx = \frac{1}{\pi} \frac{2n\sin(2n\pi x)\sin \pi x + \cos(2n\pi x)\cos \pi x}{4n^2 - 1}$$

we thereby obtain for $0 < u < 1$

$$(2.12) \quad \sum_{n=1}^{\infty} \frac{\cos(2n+1)\pi u}{2n+1} \cdot \log\left(1+\frac{1}{n}\right) = 2\sum_{n=1}^{\infty} \frac{[2n\sin(2n\pi u)\sin \pi u + \cos(2n\pi u)\cos \pi u]\log n}{4n^2-1}$$

Letting $u=0$ (or $u=1$) gives us

$$\sum_{n=1}^{\infty} \frac{1}{2n+1} \log\left(1+\frac{1}{n}\right) = 2\sum_{n=1}^{\infty} \frac{\log n}{4n^2-1}$$

and this is easily verified as follows.

We have

$$2\sum_{n=1}^{\infty} \frac{\log n}{4n^2-1} = \sum_{n=1}^{\infty} \left[\frac{1}{2n-1} - \frac{1}{2n+1}\right]\log n$$

and consider the finite sum

$$\sum_{n=1}^{N} \frac{\log n}{2n-1} = \sum_{m=0}^{N-1} \frac{\log(m+1)}{2m+1}$$

$$= \sum_{n=1}^{N} \frac{\log(n+1)}{2n+1} - \frac{\log(N+1)}{2N+1}$$

Hence we have

$$\sum_{n=1}^{N} \left[\frac{1}{2n-1} - \frac{1}{2n+1}\right]\log n = \sum_{n=1}^{N} \left[\frac{\log(n+1)}{2n+1} - \frac{\log n}{2n+1}\right] - \frac{\log(N+1)}{2N+1}$$

Therefore, as $N \to \infty$ we see that

$$2\sum_{n=1}^{\infty} \frac{\log n}{4n^2-1} = \sum_{n=1}^{\infty} \frac{1}{2n+1}\log\frac{n+1}{n}$$

$\square$

Letting $u = \frac{1}{4}$ in (2.12) and using

$$\sum_{n=1}^{\infty} a_n \cos(n\pi/2) = \sum_{n=1}^{\infty} a_{2n} \cos(n\pi)$$

$$\sum_{n=1}^{\infty} a_n \sin(n\pi/2) = \sum_{n=0}^{\infty} a_{2n+1} \cos(n\pi)$$

gives us

$$(2.13) \quad \sum_{n=1}^{\infty} \frac{(-1)^n}{4n+1} \cdot \log\left(1+\frac{1}{2n}\right) - \sum_{n=0}^{\infty} \frac{(-1)^n}{4n+3} \cdot \log\left(1+\frac{1}{2n+1}\right)$$

$$= 2\sum_{n=0}^{\infty} \frac{2(2n+1)(-1)^n \log(2n+1)}{4(2n+1)^2 - 1} + 2\sum_{n=1}^{\infty} \frac{(-1)^n \log(2n)}{16n^2 - 1}$$

□

Integrating (2.9) readily results in

(2.14) $$\int_0^1 \psi(x) \sin \pi x \, dx = -\frac{2}{\pi}\left[\gamma + \log(2\pi) + 2\sum_{n=1}^{\infty} \frac{\log n}{4n^2 - 1}\right]$$

as previously determined by Kölbig [36]. Another derivation was given in [21].

□

We may multiply (2.6) by $\sin(2k\pi x)$ and integrate to obtain a very concise derivation of the known Fourier coefficient

(2.15) $$\int_0^1 \psi(x) \sin 2k\pi x \, dx = -\frac{\pi}{2}$$

**Remark:**

Care must be exercised when we employ these 'limit' identities. For example, if we multiply (2.7) by $x$ it is not at all obvious how one could have

$$\lim_{x \to 0} \lim_{s \to 1} \sum_{n=1}^{\infty} \frac{x \sin 2n\pi x}{n^{1-s}} = \lim_{x \to 0} \frac{1}{2} x \cot(\pi x) = \frac{1}{2\pi}$$

### 3. A novel derivation of the duplication formula for the digamma function

Landau [37, p.181] showed that for $0 < x < \frac{1}{2}$

(3.1) $$f\left(x + \frac{1}{2}\right) = f(2x) - f(x) - \log 2 \log(2 \sin 2\pi x)$$

where $f(x)$ is defined as

$$f(x) := \sum_{n=1}^{\infty} \frac{\log n}{n} \cos 2n\pi x \, .$$

We may generalise Landau's formula by defining $f_s(x)$ as $f_s(x) := \sum_{n=1}^{\infty} \frac{\log n}{n^s} \cos 2n\pi x$ which gives us

$$f_s\left(x + \frac{1}{2}\right) = \sum_{n=1}^{\infty} \frac{\log n}{n^s} \cos(2n\pi x + n\pi)$$

$$= \sum_{n=1}^{\infty}(-1)^n \frac{\log n}{n^s}\cos 2n\pi x$$

Subject to the necessary convergence conditions being satisfied we have

(3.1.1) $$\sum_{n=1}^{\infty}(-1)^n a_n = -\sum_{n=1}^{\infty} a_n + 2\sum_{n=1}^{\infty} a_{2n}$$

and thus we have

$$\sum_{n=1}^{\infty}(-1)^n \frac{\log n}{n^s}\cos 2n\pi x = -\sum_{n=1}^{\infty}\frac{\log n}{n^s}\cos 2n\pi x + 2^{1-s}\sum_{n=1}^{\infty}\frac{\log(2n)}{n^s}\cos 4n\pi x$$

$$= -f_s(x) + 2^{1-s}\log 2 \sum_{n=1}^{\infty}\frac{\cos 4n\pi x}{n^s} + 2^{1-s}\sum_{n=1}^{\infty}\frac{\log n}{n^s}\cos 4n\pi x$$

Hence, we obtain

(3.2) $$f_s\left(x+\frac{1}{2}\right) = 2^{1-s} f_s(2x) - f_s(x) + 2^{1-s}\log 2 \sum_{n=1}^{\infty}\frac{\cos 4n\pi x}{n^s}$$

We easily see that

(3.3) $$\sum_{n=1}^{\infty}\frac{\cos 2n\pi(x+\frac{1}{2})}{n^s} = 2^{1-s}\sum_{n=1}^{\infty}\frac{\cos 4n\pi x}{n^s} - \sum_{n=1}^{\infty}\frac{\cos 2n\pi x}{n^s}$$

and

(3.4) $$\sum_{n=1}^{\infty}\frac{\sin 2n\pi(x+\frac{1}{2})}{n^s} = 2^{1-s}\sum_{n=1}^{\infty}\frac{\sin 4n\pi x}{n^s} - \sum_{n=1}^{\infty}\frac{\sin 2n\pi x}{n^s}$$

With regard to (2.6) we define $\psi(x,s)$ as

(3.5) $$\psi(x,s) := 2f_s(x) + 2[\gamma + \log(2\pi)]\sum_{n=1}^{\infty}\frac{\cos 2n\pi x}{n^s} - \pi\sum_{n=1}^{\infty}\frac{\sin 2n\pi x}{n^s}$$

so that $\psi(x) = \lim_{s \to 0}\psi(x,s)$.

Letting $x \to x + \frac{1}{2}$ in (3.5) results in

$$\psi\left(x+\frac{1}{2},s\right) = 2\left[2^{1-s} f_s(2x) - f_s(x) + 2^{1-s}\log 2\sum_{n=1}^{\infty}\frac{\cos 4n\pi x}{n^s}\right]$$

$$+ 2[\gamma + \log(2\pi)]\left[2^{1-s}\sum_{n=1}^{\infty}\frac{\cos 4n\pi x}{n^s} - \sum_{n=1}^{\infty}\frac{\cos 2n\pi x}{n^s}\right]$$

$$-\pi\left[2^{1-s}\sum_{n=1}^{\infty}\frac{\sin 4n\pi x}{n^s}-\sum_{n=1}^{\infty}\frac{\sin 2n\pi x}{n^s}\right]$$

Accordingly, we see that

$$\psi\left(x+\frac{1}{2},s\right)=2^{1-s}\psi(2x,s)-\psi(x,s)+2^{2-s}\log 2\sum_{n=1}^{\infty}\frac{\cos 4n\pi x}{n^s}$$

and, since $\lim_{s\to 0}\sum_{n=1}^{\infty}\frac{\cos 4n\pi x}{n^s}=-\frac{1}{2}$ in the interval $(0,\frac{1}{2})$, taking the limit as $s\to 0$ we obtain the well-known duplication formula for the digamma function [41, p.213]

(3.6) $$\psi\left(x+\frac{1}{2}\right)=2\psi(2x)-\psi(x)-2\log 2$$

## 4. A functional equation for the first generalised Stieltjes constant $\gamma_1(x)$

Following on from Section 3, in a similar way we obtain for $g_s(x):=\sum_{n=1}^{\infty}\frac{\log n}{n^s}\sin 2n\pi x$

(4.1) $$g_s\left(x+\frac{1}{2}\right)=2^{1-s}g_s(2x)-g_s(x)+2^{1-s}\log 2\sum_{n=1}^{\infty}\frac{\sin 4n\pi x}{n^s}$$

In particular, for $s=1$ we have using (5.3)

$$g_1\left(x+\frac{1}{2}\right)=-g_1(x)-\pi\left(2x-\frac{1}{2}\right)\log 2+g_1(2x)$$

which was also employed by Landau [37].

**Theorem 4.1**

We have the functional equation for the first generalised Stieltjes constant which is valid for $0<x<\frac{1}{2}$

(4.2) $$\gamma_1\left(\frac{1}{2}-x\right)-\gamma_1\left(x+\frac{1}{2}\right)=-2[\gamma_1(2x)-\gamma_1(1-2x)]+[\gamma_1(x)-\gamma_1(1-x)]+2\pi\log 2\cot(2\pi x)$$

**Proof**

We recall (1.1)

$$\gamma_1(1-x)-\gamma_1(x)=\pi[\gamma+\log(2\pi)]\cot(\pi x)+2\pi\lim_{s\to 1}\sum_{n=1}^{\infty}\frac{\sin(2n\pi x)\log n}{n^{1-s}}$$

and express this as

$$\gamma_1(1-x) - \gamma_1(x) = 2\pi[\gamma + \log(2\pi)]\lim_{s\to 0}\sum_{n=1}^{\infty}\frac{\sin(2n\pi x)}{n^s} + 2\pi\lim_{s\to 0} g_s(x)$$

Letting $x \to x + \frac{1}{2}$ gives us

$$\gamma_1\left(\frac{1}{2}-x\right) - \gamma_1\left(x+\frac{1}{2}\right) = 2\pi[\gamma + \log(2\pi)]\lim_{s\to 0}\sum_{n=1}^{\infty}(-1)^n\frac{\sin(2n\pi x)}{n^s}$$

$$+ 2\pi\lim_{s\to 0}\left[2^{1-s} g_s(2x) - g_s(x) + 2^{1-s}\log 2\sum_{n=1}^{\infty}\frac{\sin 4n\pi x}{n^s}\right]$$

Letting $x \to x + \frac{1}{2}$ in $\lim_{s\to 0}\sum_{n=1}^{\infty}\frac{\sin(2n\pi x)}{n^s} = \frac{1}{2}\cot(\pi x)$ results in

$$\lim_{s\to 0}\sum_{n=1}^{\infty}(-1)^n\frac{\sin(2n\pi x)}{n^s} = -\frac{1}{2}\tan(\pi x)$$

Hence we have

$$\gamma_1\left(\frac{1}{2}-x\right) - \gamma_1\left(x+\frac{1}{2}\right) = -\pi[\gamma + \log(2\pi)]\tan(\pi x)$$

$$+ 2\pi\lim_{s\to 0}\left[2^{1-s} g_s(2x) - g_s(x) + 2^{1-s}\log 2\sum_{n=1}^{\infty}\frac{\sin 4n\pi x}{n^s}\right]$$

Since

$$\gamma_1(1-x) - \gamma_1(x) - 2\pi[\gamma + \log(2\pi)]\lim_{s\to 0}\sum_{n=1}^{\infty}\frac{\sin(2n\pi x)}{n^s} = 2\pi\lim_{s\to 0} g_s(x)$$

we have

$$\gamma_1\left(\frac{1}{2}-x\right) - \gamma_1\left(x+\frac{1}{2}\right) = -\pi[\gamma + \log(2\pi)]\tan(\pi x)$$

$$+ 2\left[\gamma_1(1-2x) - \gamma_1(2x) - 2\pi[\gamma + \log(2\pi)]\lim_{s\to 0}\sum_{n=1}^{\infty}\frac{\sin(4n\pi x)}{n^s}\right]$$

$$- \left[\gamma_1(1-x) - \gamma_1(x) - 2\pi[\gamma + \log(2\pi)]\lim_{s\to 0}\sum_{n=1}^{\infty}\frac{\sin(2n\pi x)}{n^s}\right]$$

$$+ 2\pi\lim_{s\to 0} 2^{1-s}\log 2\sum_{n=1}^{\infty}\frac{\sin 4n\pi x}{n^s}$$

$$= -\pi[\gamma + \log(2\pi)]\tan(\pi x)$$

$$+ 2\big[\gamma_1(1-2x) - \gamma_1(2x) - \pi[\gamma + \log(2\pi)]\cot(2\pi x)\big]$$

$$- \big[\gamma_1(1-x) - \gamma_1(x) - \pi[\gamma + \log(2\pi)]\cot(\pi x)\big]$$

$$+ 2\pi \lim_{s \to 0} 2^{1-s} \log 2 \sum_{n=1}^{\infty} \frac{\sin 4n\pi x}{n^s}$$

Using the trigonometric identity $\cot(\pi x) - \tan(\pi x) = 2\cot(2\pi x)$ this simplifies to

$$\gamma_1\left(\frac{1}{2} - x\right) - \gamma_1\left(x + \frac{1}{2}\right) = -2\big[\gamma_1(2x) - \gamma_1(1-2x)\big] + \big[\gamma_1(x) - \gamma_1(1-x)\big] + 2\pi \log 2 \cot(2\pi x)$$

We obtain a more explicit result in (7.19) below.

### 5. A formal derivation of Kummer's Fourier Series for $\log \Gamma(x)$

We recall (2.6)

$$\psi(x) = \lim_{s \to 1} \sum_{n=1}^{\infty} \frac{2[\gamma + \log(2\pi n)]\cos(2n\pi x) - \pi \sin(2n\pi x)}{n^{1-s}}$$

Assuming that the following operations are valid, formal integration of this over $[\frac{1}{2}, u]$ results in (where $0 < u < 1$)

$$\int_{\frac{1}{2}}^{u} \lim_{s \to 1} \sum_{n=1}^{\infty} \frac{[\gamma + \log(2\pi n)]\cos 2n\pi x}{n^{1-s}} dx = \lim_{s \to 1} \int_{\frac{1}{2}}^{u} \sum_{n=1}^{\infty} \frac{[\gamma + \log(2\pi n)]\cos 2n\pi x}{n^{1-s}} dx$$

$$= \lim_{s \to 1} \sum_{n=1}^{\infty} \frac{\gamma + \log(2\pi n)}{n^{1-s}} \int_{\frac{1}{2}}^{u} \cos 2n\pi x \, dx$$

$$= \sum_{n=1}^{\infty} \frac{\gamma + \log(2\pi n)}{2\pi n} \sin 2n\pi u$$

In a similar fashion we have

$$\int_{\frac{1}{2}}^{u} \lim_{s \to 1} \sum_{n=1}^{\infty} \frac{\sin 2n\pi x}{n^{1-s}} dx = \lim_{s \to 1} \int_{\frac{1}{2}}^{u} \sum_{n=1}^{\infty} \frac{\sin 2n\pi x}{n^{1-s}} dx$$

$$= \lim_{s \to 1} \sum_{n=1}^{\infty} \frac{1}{n^{1-s}} \int_{\frac{1}{2}}^{u} \sin 2n\pi x \, dx$$

$$= \sum_{n=1}^{\infty} \frac{(-1)^n - \cos 2n\pi u}{2\pi n}$$

We then obtain

$$\log \Gamma(u) - \log \Gamma(\tfrac{1}{2}) = \sum_{n=1}^{\infty} \frac{\gamma + \log(2\pi n)}{\pi n} \sin 2n\pi u - \frac{1}{2} \sum_{n=1}^{\infty} \frac{(-1)^n - \cos 2n\pi u}{n}$$

and thus we have for $0 < u < 1$

(5.1) $$\log \Gamma(u) = \frac{1}{2} \log \pi + \sum_{n=1}^{\infty} \frac{\gamma + \log(2\pi n)}{\pi n} \sin 2n\pi u + \frac{1}{2} \sum_{n=1}^{\infty} \frac{\cos 2n\pi u}{n} + \frac{1}{2} \log 2$$

Using the familiar Fourier series shown in Carslaw's book [15, p.241]

(5.2) $$\log[2 \sin \pi u] = -\sum_{n=1}^{\infty} \frac{\cos 2n\pi u}{n} \qquad (0 < u < 1)$$

(5.3) $$\pi \left( u - \frac{1}{2} \right) = -\sum_{n=1}^{\infty} \frac{\sin 2n\pi u}{n} \qquad (0 < u < 1)$$

we may express (5.1) as Kummer's formula [44]

(5.4) $$\log \Gamma(u) = \frac{1}{2} \log \pi + \sum_{n=1}^{\infty} \frac{\log n}{\pi n} \sin 2n\pi u - \left( u - \frac{1}{2} \right)[\gamma + \log(2\pi)] - \frac{1}{2} \log \sin \pi u$$

or equivalently

(5.5) $$\log \Gamma(u) = \frac{1}{2} \log \pi + \sum_{n=1}^{\infty} \frac{[\gamma + \log(2\pi n)] \sin 2\pi nu}{\pi n} - \frac{1}{2} \log \sin \pi u$$

It may be noted that one may also formally derive (5.2) and (5.3) by integrating (1.4) and (1.5) over the interval $[\tfrac{1}{2}, u]$.

Alternatively, let us differentiate (2.1) to obtain

$$\varsigma'(s,x) = 2\sum_{n=1}^{\infty} \left[ \log(2\pi n) - \psi(1-s) + \frac{\pi}{2} \cot\left(\frac{\pi s}{2}\right) \right] \frac{\Gamma(1-s)}{(2\pi n)^{1-s}} \sin\left(\frac{\pi s}{2}\right) \cos 2n\pi x$$

$$+ 2\sum_{n=1}^{\infty} \left[ \log(2\pi n) - \psi(1-s) - \frac{\pi}{2} \tan\left(\frac{\pi s}{2}\right) \right] \frac{\Gamma(1-s)}{(2\pi n)^{1-s}} \cos\left(\frac{\pi s}{2}\right) \sin 2n\pi x$$

*Assuming* this is valid at $s = 0$ we have

$$\varsigma'(0,x) = \frac{1}{2}\sum_{n=1}^{\infty}\frac{\cos 2n\pi x}{n} + \frac{1}{\pi}\sum_{n=1}^{\infty}\frac{[\gamma+\log(2\pi n)]}{n}\sin 2\pi nx$$

Then using Lerch's identity [10] for $x > 0$

(5.6) $$\varsigma'(0,x) = \log\Gamma(x) - \frac{1}{2}\log(2\pi)$$

we obtain another derivation of Kummer's formula.

This was the method applied by Berndt [10] who has indeed pointed out that Hurwitz's formula for the Fourier expansion of the Riemann zeta function $\varsigma(s,x)$ also holds for $\text{Re}(s) < 1$ with $0 < x < 1$.

## 6. A representation involving the Barnes double gamma function $G(x)$

**Theorem**

We have for $0 < u < 1$

(6.1) $$\log G(1+u) = \frac{1}{2}u\log(2\pi) - \frac{1}{2}u(u-1) - \frac{1}{4\pi}\sum_{n=1}^{\infty}\frac{\sin 2n\pi u}{n^2} + \frac{1}{2}u\sum_{n=1}^{\infty}\frac{\cos 2n\pi u}{n}$$

$$+u\sum_{n=1}^{\infty}\frac{[\gamma+\log(2\pi n)]\sin(2n\pi u)}{n\pi} + \frac{1}{2\pi^2}\sum_{n=1}^{\infty}\frac{[\gamma+\log(2\pi n)][\cos(2n\pi u)-1]}{n^2}$$

where $G(x)$, the Barnes double gamma function, is defined by [44, p.25]

(6.2) $$\log G(1+x) = \frac{1}{2}x\log(2\pi) - \frac{1}{2}x(1+x) - \frac{1}{2}\gamma x^2 + \sum_{n=1}^{\infty}\left[\frac{1}{2n}x^2 - x + n\log\left(1+\frac{x}{n}\right)\right]$$

**Proof**

We multiply (2.6) by $x$ and integrate to obtain

$$\int_0^u x\psi(x)dx = -\frac{\pi}{2}\int_0^u x\cot(\pi x)dx + \lim_{s\to 1}2\int_0^u x\sum_{n=1}^{\infty}\frac{[\gamma+\log(2\pi n)]\cos(2n\pi x)}{n^{1-s}}dx$$

We have

$$\int_0^u x\cos(2n\pi x)dx = \frac{x\sin(2n\pi x)}{(2n\pi)} + \frac{\cos(2n\pi x)}{(2n\pi)^2}\bigg|_0^u$$

$$= \frac{u \sin(2n\pi u)}{(2n\pi)} + \frac{\cos(2n\pi u) - 1}{(2n\pi)^2}$$

which gives us

$$\lim_{s \to 1} 2 \int_0^u x \sum_{n=1}^{\infty} \frac{[\gamma + \log(2\pi n)] \cos(2n\pi x)}{n^{1-s}} dx = 2u \sum_{n=1}^{\infty} \frac{[\gamma + \log(2\pi n)] \sin(2n\pi u)}{(2n\pi)}$$

$$+ 2 \sum_{n=1}^{\infty} \frac{[\gamma + \log(2\pi n)][\cos(2n\pi u) - 1]}{(2n\pi)^2}$$

An elementary proof of the following integral was provided in [25]

(6.3) $$\int_0^u x \psi(1+x) dx = \log G(1+u) - \frac{1}{2} u \log(2\pi) + \frac{1}{2} u(u+1)$$

and, since $\psi(1+x) = \psi(x) + \frac{1}{x}$, we therefore have

(6.4) $$\int_0^u x \psi(x) dx = \log G(1+u) - \frac{1}{2} u \log(2\pi) + \frac{1}{2} u(u-1)$$

In [18] we proved the basic identity

(6.5) $$\int_a^b p(x) \cot(\alpha x / 2) dx = 2 \sum_{n=1}^{\infty} \int_a^b p(x) \sin \alpha nx \, dx$$

where we initially required that $p(x)$ is a twice continuously differentiable function. Equation (6.5) is valid provided (i) $\sin(\alpha x/2) \neq 0 \ \forall \ x \in [a,b]$ or, alternatively, (ii) if $\sin(\alpha \eta / 2) = 0$ for some $\eta \in [a,b]$ then $p(\eta) = 0$.

As later shown in [22], (6.5) is actually valid for a wider class of suitably behaved functions (which are not necessarily twice continuously differentiable).

With $p(x) = x$ and $\alpha = 2\pi$ in (6.5) we have for $-1 < u < 1$

$$\int_0^u x \cot \pi x \, dx = 2 \sum_{n=1}^{\infty} \int_0^u x \sin 2\pi nx \, dx$$

and since

$$\int_0^u x \sin 2\pi nx \, dx = \frac{\sin 2\pi nu}{(2\pi n)^2} - \frac{u \cos 2\pi nu}{2\pi n}$$

we obtain

(6.6) $$\int_0^u x \cot \pi x \, dx = \frac{1}{2\pi^2} \sum_{n=1}^{\infty} \frac{\sin 2n\pi u}{n^2} - \frac{u}{\pi} \sum_{n=1}^{\infty} \frac{\cos 2n\pi u}{n}$$

Therefore we obtain (6.1)

$$\log G(1+u) = \frac{1}{2}u\log(2\pi) - \frac{1}{2}u(u-1) - \frac{1}{4\pi}\sum_{n=1}^{\infty}\frac{\sin 2n\pi u}{n^2} + \frac{1}{2}u\sum_{n=1}^{\infty}\frac{\cos 2n\pi u}{n}$$

$$+ u\sum_{n=1}^{\infty}\frac{[\gamma + \log(2\pi n)]\sin(2n\pi u)}{n\pi} + \frac{1}{2\pi^2}\sum_{n=1}^{\infty}\frac{[\gamma + \log(2\pi n)][\cos(2n\pi u) - 1]}{n^2}$$

$\square$

Letting $u = \frac{1}{2}$ in (6.1) we have

$$\log G(1+\tfrac{1}{2}) = \frac{1}{4}\log(2\pi) + \frac{1}{8} - \frac{1}{4}\sum_{n=1}^{\infty}\frac{(-1)^{n+1}}{n} + \frac{1}{2\pi^2}\sum_{n=1}^{\infty}(-1)^n\frac{\gamma + \log(2\pi n)}{n^2} - \frac{1}{2\pi^2}\sum_{n=1}^{\infty}\frac{\gamma + \log(2\pi n)}{n^2}$$

$$= \frac{1}{4}\log(2\pi) + \frac{1}{8} - \frac{1}{4}\log 2 + \frac{1}{2\pi^2}[\varsigma_a'(2) + \varsigma'(2)] - \frac{1}{2\pi^2}[\gamma + \log(2\pi)][\varsigma_a(2) + \varsigma(2)]$$

where $\varsigma_a(s) = \sum_{n=1}^{\infty}\frac{(-1)^{n+1}}{n^s}$ is the alternating Riemann zeta function.

We have

$$\varsigma_a(s) = (1 - 2^{1-s})\varsigma(s)$$

and hence

$$\varsigma_a'(s) = (1 - 2^{1-s})\varsigma'(s) + 2^{1-s}\log 2 . \varsigma(s)$$

We see that

$$\varsigma_a'(2) = \frac{1}{2}[\varsigma'(2) + \log 2 . \varsigma(2)]$$

We substitute the well-known result (which follows from the functional equation for the Riemann zeta function)

$$\varsigma'(-1) - \frac{1}{12}(1 - \gamma - \log 2\pi) = \frac{1}{2\pi^2}\varsigma'(2)$$

and employing [44, p.25]

$$G(1+u) = G(u)\Gamma(u)$$

we obtain

(6.7) $$\log G(\tfrac{1}{2}) = \frac{1}{24}\log 2 - \frac{1}{4}\log \pi + \frac{3}{2}\varsigma'(-1)$$

which was originally derived by Barnes [44, p.26].

□

Since

$$B_{2N}(t) = (-1)^{N+1} 2(2N)! \sum_{n=1}^{\infty} \frac{\cos 2n\pi t}{(2\pi n)^{2N}} \qquad , N = 1, 2, \ldots$$

we have

$$B_2(t) = \frac{1}{\pi^2} \sum_{n=1}^{\infty} \frac{\cos 2n\pi t}{n^2}$$

It is well known that

(6.8) $$B_2(u) = u^2 - u + \frac{1}{6}$$

Substituting (5.3) in (6.1) and using (6.8) gives us

(6.9) $$\log G(1+t) - t \log \Gamma(t) = -\frac{1}{4\pi} \sum_{n=1}^{\infty} \frac{\sin 2n\pi t}{n^2} + \frac{1}{2\pi^2} \left( \log(2\pi) + \gamma - 1 \right) \sum_{n=1}^{\infty} \frac{\cos 2n\pi t}{n^2}$$

$$+ \frac{1}{2\pi^2} \sum_{n=1}^{\infty} \frac{\log n}{n^2} \cos 2n\pi t + \varsigma'(-1)$$

We therefore have the Fourier series which we obtained by a more orthodox method in [17]. It may be noted that Mező [40] also obtained this Fourier series by an entirely different method.

We recall that the Clausen function $\text{Cl}_2(u)$ is defined by

$$\text{Cl}_2(u) = \sum_{n=1}^{\infty} \frac{\sin nu}{n^2} = -\int_0^u \log\left( 2 \sin \frac{x}{2} \right) dx$$

By differentiating (2.1) we obtain with $s = -1$

(6.10) $$\varsigma'(-1,t) = \frac{1}{4\pi} \sum_{n=1}^{\infty} \frac{\sin 2n\pi t}{n^2} - \frac{1}{2\pi^2} \left[ \log(2\pi) + \gamma - 1 \right] \sum_{n=1}^{\infty} \frac{\cos 2n\pi t}{n^2} - \frac{1}{2\pi^2} \sum_{n=1}^{\infty} \frac{\log n}{n^2} \cos 2n\pi t$$

Comparing (6.9) with (6.10) we see that

(6.11) $$\log G(1+t) - t \log \Gamma(t) = \varsigma'(-1) - \varsigma'(-1,t)$$

This functional equation was derived by Vardi [48] in 1988 and also by Gosper [31] in 1997 (an elementary derivation is contained in [25]).

## 7. A new representation for the first Stieltjes constant

**Theorem**

We have for $0 < x < 1$

(7.1) $$\gamma_1(x) = \lim_{s \to 0} \gamma_1(x, s)$$

where $\gamma_1(x, s)$ is defined by

(7.2) $$\gamma_1(x, s) = \sum_{n=1}^{\infty} \frac{[\gamma + \log(2\pi n)]^2 \cos(2\pi nx)}{n^s} - \frac{1}{2}\varsigma(2)\sum_{n=1}^{\infty} \frac{\cos(2\pi nx)}{n^s}$$

$$-\pi \sum_{n=1}^{\infty} \frac{[\gamma + \log(2\pi n)]\sin(2\pi nx)}{n^s}$$

**Proof**

Differentiating (2.5) gives us

(7.3) $$\frac{\partial^2}{\partial s^2}(s-1)[\varsigma(s, x)]$$

$$= -2\Gamma(2-s)\sum_{n=1}^{\infty} \frac{\log^2(2\pi n)\sin\left(2n\pi x + \frac{\pi s}{2}\right) + \pi \log(2\pi n)\cos\left(2n\pi x + \frac{\pi s}{2}\right) - \left(\frac{\pi}{2}\right)^2 \sin\left(2n\pi x + \frac{\pi s}{2}\right)}{(2\pi n)^{1-s}}$$

$$+ 4\Gamma'(2-s)\sum_{n=1}^{\infty} \frac{\log(2\pi n)\sin\left(2n\pi x + \frac{\pi s}{2}\right) + \frac{\pi}{2}\cos\left(2n\pi x + \frac{\pi s}{2}\right)}{(2\pi n)^{1-s}}$$

$$- 2\Gamma''(2-s)\sum_{n=1}^{\infty} \frac{\sin\left(2n\pi x + \frac{\pi s}{2}\right)}{(2\pi n)^{1-s}}$$

With $s = 1$ we have

$$\left.\frac{\partial^2}{\partial s^2}(s-1)[\varsigma(s, x)]\right|_{s=1}$$

$$= -2\lim_{s \to 1}\sum_{n=1}^{\infty} \frac{\log^2(2\pi n)\cos(2n\pi x) - \pi \log(2\pi n)\sin(2n\pi x) - \left(\frac{\pi}{2}\right)^2 \cos(2n\pi x)}{n^{1-s}}$$

$$-4\gamma \lim_{s \to 1}\sum_{n=1}^{\infty} \frac{\log(2\pi n)\cos(2n\pi x) - \frac{\pi}{2}\sin(2n\pi x)}{n^{1-s}} - 2[\gamma^2 + \varsigma(2)]\lim_{s \to 1}\sum_{n=1}^{\infty} \frac{\cos(2n\pi x)}{n^{1-s}}$$

since $\Gamma''(1) = \gamma^2 + \varsigma(2)$.

We note that

$$\frac{\partial^2}{\partial s^2}(s-1)[\varsigma(s,x)]\bigg|_{s=1} = -2\gamma_1(x)$$

and we accordingly define $\gamma_1(x,s)$ as

$$\gamma_1(x,s) := \sum_{n=1}^{\infty} \frac{\log^2(2\pi n)\cos(2n\pi x) - \pi\log(2\pi n)\sin(2n\pi x) - \left(\frac{\pi}{2}\right)^2\cos(2n\pi x)}{n^s}$$

$$+ 2\gamma\sum_{n=1}^{\infty} \frac{\log(2\pi n)\cos(2n\pi x) - \frac{\pi}{2}\sin(2n\pi x)}{n^s}$$

$$+ [\gamma^2 + \varsigma(2)]\sum_{n=1}^{\infty} \frac{\cos(2n\pi x)}{n^s}$$

and we see that $\lim_{s \to 0} \gamma_1(x,s) = \gamma_1(x)$.

With a little algebra we may express this as

(7.4) $\quad \gamma_1(x,s) = \sum_{n=1}^{\infty} \frac{[\gamma + \log(2\pi n)]^2 \cos(2\pi nx)}{n^s} - \frac{1}{2}\varsigma(2)\sum_{n=1}^{\infty} \frac{\cos(2\pi nx)}{n^s}$

$$-\pi\sum_{n=1}^{\infty} \frac{[\gamma + \log(2\pi n)]\sin(2\pi nx)}{n^s}$$

This may be generalised to determine $\gamma_n(x,s)$ by taking higher derivatives of (7.3).

**Determination of $\gamma_1$**

We have with $x = \frac{1}{2}$

$$\gamma_1(\tfrac{1}{2},s) = \sum_{n=1}^{\infty} \frac{(-1)^n \log^2 n}{n^s} + 2[\gamma + \log(2\pi)]\sum_{n=1}^{\infty} \frac{(-1)^n \log n}{n^s}$$

$$+ \left[[\gamma + \log(2\pi)]^2 - \frac{1}{2}\varsigma(2)\right]\sum_{n=1}^{\infty} \frac{(-1)^n}{n^s}$$

and in the limit $s \to 0$ we have

$$\gamma_1(\tfrac{1}{2}) = -\varsigma_a''(0) + 2[\gamma + \log(2\pi)]\varsigma_a'(0) - \left[[\gamma + \log(2\pi)]^2 - \tfrac{1}{2}\varsigma(2)\right]\varsigma_a(0)$$

Substituting

$$\varsigma_a(0) = \tfrac{1}{2}$$

$$\varsigma_a'(0) = \frac{1}{2}\log\left(\frac{\pi}{2}\right)$$

$$\varsigma_a''(0) = -\varsigma''(0) - 2\log 2 \log(2\pi) + \log^2 2$$

we obtain

$$\gamma_1(\tfrac{1}{2}) = \varsigma''(0) + 2\log 2 \log(2\pi) - \log^2 2$$

$$+ [\gamma + \log(2\pi)]\log\left(\frac{\pi}{2}\right) - \frac{1}{2}\left[[\gamma + \log(2\pi)]^2 - \frac{1}{2}\varsigma(2)\right]$$

It is known that for integers $q \geq 2$

(7.5) $$\sum_{r=1}^{q-1}\gamma_p\left(\frac{r}{q}\right) = -\gamma_p + q(-1)^p \frac{\log^{p+1} q}{p+1} + q\sum_{j=0}^{p}\binom{p}{j}(-1)^j \gamma_{p-j}\log^j q$$

which was previously derived by Coffey [17]. The relationship (7.5) has some history in that the case $p = 1$ was originally discovered by Ramanujan [7, p.198].

Hence we obtain

(7.6) $$\gamma_1(\tfrac{1}{2}) = \gamma_1 - \log^2 2 - 2\gamma \log 2$$

Using (7.6), a little algebra gives us

(7.7) $$\gamma_1 = \varsigma''(0) + [\gamma + \log(2\pi)]\log(2\pi) - \frac{1}{2}\left([\gamma + \log(2\pi)]^2 - \frac{1}{2}\varsigma(2)\right)$$

which corresponds with the result previously obtained by Ramanujan [9] and, at a later date, by Apostol [5] in the equivalent form

$$\varsigma''(0) = \gamma_1 + \frac{1}{2}\gamma^2 - \frac{1}{24}\pi^2 - \frac{1}{2}\log^2(2\pi)$$

**Determination of** $\gamma_1\left(\tfrac{1}{4}\right)$

Writing (7.4) as

$$(7.8) \qquad \gamma_1(x,s) = \sum_{n=1}^{\infty} \frac{\log^2 n \cos(2n\pi x)}{n^s} + 2[\gamma + \log(2\pi)] \sum_{n=1}^{\infty} \frac{\log n \cos(2n\pi x)}{n^s}$$

$$+ \left[ [\gamma + \log(2\pi)]^2 - \frac{1}{2}\varsigma(2) \right] \sum_{n=1}^{\infty} \frac{\cos(2n\pi x)}{n^s}$$

$$-\pi \sum_{n=1}^{\infty} \frac{\log n \sin(2n\pi x)}{n^s} - \pi[\gamma + \log(2\pi)] \sum_{n=1}^{\infty} \frac{\sin(2n\pi x)}{n^s}$$

we see that

$$\gamma_1(\tfrac{1}{4},s) = \sum_{n=1}^{\infty} \frac{\log^2 n \cos(n\pi/2)}{n^s} + 2[\gamma + \log(2\pi)] \sum_{n=1}^{\infty} \frac{\log n \cos(n\pi/2)}{n^s}$$

$$+ \left[ [\gamma + \log(2\pi)]^2 - \frac{1}{2}\varsigma(2) \right] \sum_{n=1}^{\infty} \frac{\cos(n\pi/2)}{n^s}$$

$$-\pi \sum_{n=1}^{\infty} \frac{\log n \sin(n\pi/2)}{n^s} - \pi[\gamma + \log(2\pi)] \sum_{n=1}^{\infty} \frac{\sin(n\pi/2)}{n^s}$$

We have for suitably convergent series

$$\sum_{n=1}^{\infty} a_n \cos(n\pi/2) = \sum_{n=1}^{\infty} a_{2n} \cos(n\pi)$$

$$\sum_{n=1}^{\infty} a_n \sin(n\pi/2) = \sum_{n=0}^{\infty} a_{2n+1} \cos(n\pi)$$

and thus

$$\gamma_1(\tfrac{1}{4},s) = \sum_{n=1}^{\infty} \frac{\log^2(2n) \cos(n\pi)}{(2n)^s} + 2[\gamma + \log(2\pi)] \sum_{n=1}^{\infty} \frac{\log(2n) \cos(n\pi)}{(2n)^s}$$

$$+ \left[ [\gamma + \log(2\pi)]^2 - \frac{1}{2}\varsigma(2) \right] \sum_{n=1}^{\infty} \frac{\cos(n\pi)}{(2n)^s}$$

$$-\pi \sum_{n=0}^{\infty} \frac{\log(2n+1) \cos(n\pi)}{(2n+1)^s} - \pi[\gamma + \log(2\pi)] \sum_{n=0}^{\infty} \frac{\cos(n\pi)}{(2n+1)^s}$$

This may be expressed as

$$\gamma_1(\tfrac{1}{4},s) = \frac{1}{2^s} \sum_{n=1}^{\infty} \frac{(-1)^n \log^2(n)}{n^s} + \frac{1}{2^s} 2[\gamma + \log(2\pi) + \log 2] \sum_{n=1}^{\infty} \frac{(-1)^n \log(n)}{n^s}$$

$$+\left[[\gamma+\log(2\pi)]^2 - \frac{1}{2}\varsigma(2) + \log^2 2 + 2[\gamma+\log(2\pi)]\log 2\right]\frac{1}{2^s}\sum_{n=1}^{\infty}\frac{(-1)^n}{n^s}$$

$$-\pi\sum_{n=0}^{\infty}\frac{(-1)^n \log(2n+1)}{(2n+1)^s} - \pi[\gamma+\log(2\pi)]\sum_{n=0}^{\infty}\frac{(-1)^n}{(2n+1)^s}$$

Employing the alternating Riemann zeta function this may equivalently be written as

$$\gamma_1(\tfrac{1}{4}, s) = -\frac{1}{2^s}\varsigma_a''(s) + \frac{1}{2^s}2[\gamma+\log(2\pi)+\log 2]\varsigma_a'(s)$$

$$-\left[[\gamma+\log(2\pi)]^2 - \frac{1}{2}\varsigma(2) + \log^2 2 + 2[\gamma+\log(2\pi)]\log 2\right]\frac{1}{2^s}\varsigma_a(s)$$

$$-\pi\sum_{n=0}^{\infty}\frac{(-1)^n \log(2n+1)}{(2n+1)^s} - \pi[\gamma+\log(2\pi)]\sum_{n=0}^{\infty}\frac{(-1)^n}{(2n+1)^s}$$

Therefore we obtain

$$\gamma_1(\tfrac{1}{4}) = -\varsigma_a''(0) + 2[\gamma+\log(2\pi)+\log 2]\varsigma_a'(0)$$

$$-\left[[\gamma+\log(2\pi)]^2 - \frac{1}{2}\varsigma(2) + \log^2 2 + 2[\gamma+\log(2\pi)]\log 2\right]\varsigma_a(0)$$

$$-\pi\lim_{s\to 0}\sum_{n=0}^{\infty}\frac{(-1)^n \log(2n+1)}{(2n+1)^s} - \pi[\gamma+\log(2\pi)]\lim_{s\to 0}\sum_{n=0}^{\infty}\frac{(-1)^n}{(2n+1)^s}$$

which becomes

$$\gamma_1(\tfrac{1}{4}) = \varsigma''(0) + 2\log 2\log(2\pi) - \log^2 2 + [\gamma+\log(2\pi)+\log 2]\log\left(\frac{\pi}{2}\right)$$

$$-\frac{1}{2}\left[[\gamma+\log(2\pi)]^2 - \frac{1}{2}\varsigma(2) + \log^2 2 + 2[\gamma+\log(2\pi)]\log 2\right]$$

$$-\pi\lim_{s\to 0}\sum_{n=0}^{\infty}\frac{(-1)^n \log(2n+1)}{(2n+1)^s} - \pi[\gamma+\log(2\pi)]\lim_{s\to 0}\sum_{n=0}^{\infty}\frac{(-1)^n}{(2n+1)^s}$$

Using (7.7) we obtain

(7.9) $$\gamma_1(\tfrac{1}{4}) = \frac{1}{2}[2\gamma_1 - 7\log^2 2 - 6\gamma\log 2]$$

$$-\pi\lim_{s\to 0}\sum_{n=0}^{\infty}\frac{(-1)^n \log(2n+1)}{(2n+1)^s} - \pi[\gamma+\log(2\pi)]\lim_{s\to 0}\sum_{n=0}^{\infty}\frac{(-1)^n}{(2n+1)^s}$$

We now proceed to determine the other two limits.

Upon a separation of terms according to the parity of $n$ we see that for Re $s > 1$

$$\varsigma_a(s,t) = \sum_{n=0}^{\infty} \frac{(-1)^n}{(n+t)^s} = \sum_{n=0}^{\infty} \frac{1}{(2n+t)^s} - \sum_{n=0}^{\infty} \frac{1}{(2n+1+t)^s}$$

$$= 2^{-s}\left[\sum_{n=0}^{\infty} \frac{1}{(n+t/2)^s} - \sum_{n=0}^{\infty} \frac{1}{(n+(1+t)/2)^s}\right]$$

and we therefore see that $\varsigma_a(s,t)$ is related to the Hurwitz zeta function by the formula

(7.10) $$\varsigma_a(s,t) = 2^{-s}\left[\varsigma\left(s,\frac{t}{2}\right) - \varsigma\left(s,\frac{1+t}{2}\right)\right]$$

The Dirichlet beta function $\beta(s)$ is defined for Re$(s) > 0$ by

(7.11) $$\beta(s) = \sum_{n=0}^{\infty} \frac{(-1)^n}{(2n+1)^s}$$

or equivalently as

(7.12) $$\beta(s) = \frac{1}{2^s}\varsigma_a(s, \tfrac{1}{2})$$

Therefore, using (7.10), the following definition, in terms of the Hurwitz zeta function, is valid in the whole complex s-plane

(7.13) $$\beta(s) = \frac{1}{4^s}\left[\varsigma\left(s,\frac{1}{4}\right) - \varsigma\left(s,\frac{3}{4}\right)\right]$$

Differentiation results in

(7.13.1) $$\beta'(s) = \sum_{n=0}^{\infty} \frac{(-1)^{n+1}\log(2n+1)}{(2n+1)^s}$$

(7.13.2) $$\beta'(s) = \frac{1}{4^s}\left[\varsigma'\left(s,\frac{1}{4}\right) - \varsigma'\left(s,\frac{3}{4}\right)\right] - \frac{\log 4}{4^s}\left[\varsigma\left(s,\frac{1}{4}\right) - \varsigma\left(s,\frac{3}{4}\right)\right]$$

and hence

$$\beta'(0) = \varsigma'\left(0,\frac{1}{4}\right) - \varsigma'\left(0,\frac{3}{4}\right) - \log 4\left[\varsigma\left(0,\frac{1}{4}\right) - \varsigma\left(0,\frac{3}{4}\right)\right]$$

Since $\varsigma(0,x) = \frac{1}{2} - x$ we have

$$\beta(0) = \frac{1}{2}$$

and

$$\beta'(0) = \varsigma'\left(0, \frac{1}{4}\right) - \varsigma'\left(0, \frac{3}{4}\right) - \log 2$$

Using Lerch's identity we obtain

$$\beta'(0) = \log \Gamma\left(\frac{1}{4}\right) - \log \Gamma\left(\frac{3}{4}\right) - \log 2$$

$$= 2\log \Gamma\left(\frac{1}{4}\right) - \log \pi - \frac{3}{2}\log 2$$

Therefore, having regard to (5.9), we have

(7.14) $$\lim_{s \to 0} \sum_{n=0}^{\infty} \frac{(-1)^n}{(2n+1)^s} = -\beta(0)$$

(7.15) $$\lim_{s \to 0} \sum_{n=0}^{\infty} \frac{(-1)^n \log(2n+1)}{(2n+1)^s} = -\beta'(0)$$

and we conclude that

(7.16) $$\gamma_1\left(\frac{1}{4}\right) = \frac{1}{2}[2\gamma_1 - 7\log^2 2 - 6\gamma \log 2] - \frac{1}{2}\pi\left[\gamma + 4\log 2 + 3\log \pi - 4\log \Gamma\left(\frac{1}{4}\right)\right]$$

This is now a well-known result (see [11] and the various references therein).

We see that

(7.17) $$\gamma_1(x,s) - \gamma_1(1-x,s) = -2\pi \sum_{n=1}^{\infty} \frac{[\gamma + \log(2\pi n)]}{n^s} \sin(2n\pi x)$$

(7.18) $$\gamma_1(x,s) + \gamma_1(1-x,s) = \sum_{n=1}^{\infty} \frac{[\gamma + \log(2\pi n)]^2}{n^s} \cos(2n\pi x) - \varsigma(2)\sum_{n=1}^{\infty} \frac{\cos(2n\pi x)}{n^s}$$

which explains the 'shape' of

$$\gamma_1\left(\frac{1}{4}\right) = \frac{1}{2}[2\gamma_1 - 7\log^2 2 - 6\gamma \log 2] - \frac{1}{2}\pi\left[\gamma + 4\log 2 + 3\log \pi - 4\log \Gamma\left(\frac{1}{4}\right)\right]$$

and

$$\gamma_1\left(\frac{3}{4}\right) = \frac{1}{2}[2\gamma_1 - 7\log^2 2 - 6\gamma \log 2] + \frac{1}{2}\pi\left[\gamma + 4\log 2 + 3\log \pi - 4\log \Gamma\left(\frac{1}{4}\right)\right]$$

We see from (7.13.2) that

$$\beta'(s) = \frac{1}{4^s}\left[\varsigma'\left(s,\frac{1}{4}\right) - \varsigma'\left(s,\frac{3}{4}\right)\right] - \frac{\log 4}{4^s}\left[\varsigma\left(s,\frac{1}{4}\right) - \varsigma\left(s,\frac{3}{4}\right)\right]$$

and therefore using (1.2.2) we obtain

$$\beta'(1) = -\frac{1}{4}\left[\gamma_1\left(\frac{1}{4}\right) - \gamma_1\left(\frac{3}{4}\right)\right] + \frac{\log 4}{4}\left[\psi\left(\frac{1}{4}\right) - \psi\left(\frac{3}{4}\right)\right]$$

Hence we have

$$\beta'(1) = \frac{1}{4}\pi\left[\gamma + 2\log 2 + 3\log\pi - 4\log\Gamma\left(\frac{1}{4}\right)\right]$$

$$= \sum_{n=0}^{\infty}\frac{(-1)^{n+1}\log(2n+1)}{2n+1}$$

as originally determined by Malmstén [11] in 1842.

**Theorem**

We have the functional equation for the first Stieltjes constant

(7.19) $\quad \gamma_1(x+\tfrac{1}{2}) = 2\gamma_1(2x) - \gamma_1(x) - \log^2 2 + 2\log 2 \cdot \psi(2x)$

**Proof**

We recall (7.8)

$$\gamma_1(x,s) = \sum_{n=1}^{\infty}\frac{\log^2 n \cos(2n\pi x)}{n^s} + 2[\gamma+\log(2\pi)]\sum_{n=1}^{\infty}\frac{\log n \cos(2n\pi x)}{n^s}$$

$$+ \left[[\gamma+\log(2\pi)]^2 - \frac{1}{2}\varsigma(2)\right]\sum_{n=1}^{\infty}\frac{\cos(2n\pi x)}{n^s}$$

$$-\pi\sum_{n=1}^{\infty}\frac{\log n \sin(2n\pi x)}{n^s} - \pi[\gamma+\log(2\pi)]\sum_{n=1}^{\infty}\frac{\sin(2n\pi x)}{n^s}$$

which, for convenience, we write as

$$\gamma_1(x,s) = h_s(x) - \pi g_s(x) + 2[\gamma+\log(2\pi)]f_s(x)$$

$$-\pi[\gamma+\log(2\pi)]S_s(x) + \left[[\gamma+\log(2\pi)]^2 - \frac{1}{2}\varsigma(2)\right]C_s(x)$$

where
$$h_s(x) := \sum_{n=1}^{\infty} \frac{\log^2 n \cos(2n\pi x)}{n^s}$$

$$g_s(x) := \sum_{n=1}^{\infty} \frac{\log n \cos(2n\pi x)}{n^s}$$

$$C_s(x) := \sum_{n=1}^{\infty} \frac{\cos 2n\pi x}{n^s}$$

$$S_s(x) := \sum_{n=1}^{\infty} \frac{\sin 2n\pi x}{n^s}$$

It is easily seen that
$$C_s(x + \tfrac{1}{2}) = 2^{1-s} C_s(2x) - C_s(x)$$
$$S_s(x + \tfrac{1}{2}) = 2^{1-s} S_s(2x) - S_s(x)$$

and (3.2) and (4.1) may be written as
$$f_s\left(x + \frac{1}{2}\right) = 2^{1-s} f_s(2x) - f_s(x) + 2^{1-s} \log 2 \cdot C_s(2x)$$

$$g_s\left(x + \frac{1}{2}\right) = 2^{1-s} g_s(2x) - g_s(x) + 2^{1-s} \log 2 \cdot S_s(2x)$$

Letting $x \to x + \tfrac{1}{2}$ gives us

$$\gamma_1(x + \tfrac{1}{2}, s) = h_s(x + \tfrac{1}{2}) - \pi \left[ 2^{1-s} g_s(2x) - g_s(x) + 2^{1-s} \log 2 \cdot S_s(2x) \right]$$

$$+ 2[\gamma + \log(2\pi)] \left[ 2^{1-s} f_s(2x) - f_s(x) + 2^{1-s} \log 2 \cdot C_s(2x) \right]$$

$$- \pi[\gamma + \log(2\pi)] \left[ 2^{1-s} S_s(2x) - S_s(x) \right]$$

$$+ \left[ [\gamma + \log(2\pi)]^2 - \frac{1}{2} \varsigma(2) \right] \left[ 2^{1-s} C_s(2x) - C_s(x) \right]$$

We have
$$h_s(x + \tfrac{1}{2}) = 2^{1-s} \sum_{n=1}^{\infty} \frac{\log^2(2n) \cos(4n\pi x)}{n^s} - \sum_{n=1}^{\infty} \frac{\log^2 n \cos(2n\pi x)}{n^s}$$

$$= 2^{1-s} \sum_{n=1}^{\infty} \frac{\log^2 n \cos(4n\pi x)}{n^s} + 2^{1-s} \log^2 2 \sum_{n=1}^{\infty} \frac{\cos(4n\pi x)}{n^s} + 2^{2-s} \log 2 \sum_{n=1}^{\infty} \frac{\log n \cos(4n\pi x)}{n^s} - h_s(x)$$

$$= 2^{1-s} h_s(2x) - h_s(x) + 2^{1-s} \log^2 2 \cdot C_s(2x) + 2^{2-s} \log 2 \sum_{n=1}^{\infty} \frac{\log n \cos(4n\pi x)}{n^s}$$

and therefore we have

$$\gamma_1(x + \tfrac{1}{2}, s) = \left[ 2^{1-s} h_s(2x) - h_s(x) \right] + 2^{1-s} \log^2 2 \cdot C_s(2x) + 2^{2-s} \log 2 \sum_{n=1}^{\infty} \frac{\log n \cos(4n\pi x)}{n^s}$$

$$-\pi \left[ 2^{1-s} g_s(2x) - g_s(x) + 2^{1-s} \log 2 \cdot S_s(2x) \right]$$

$$+ 2[\gamma + \log(2\pi)] \left[ 2^{1-s} f_s(2x) - f_s(x) + 2^{1-s} \log 2 \cdot C_s(2x) \right]$$

$$-\pi[\gamma + \log(2\pi)] \left[ 2^{1-s} S_s(2x) - S_s(x) \right]$$

$$+ \left[ [\gamma + \log(2\pi)]^2 - \frac{1}{2}\varsigma(2) \right] \left[ 2^{1-s} C_s(2x) - C_s(x) \right]$$

Simple algebra results in

$$\gamma_1(x + \tfrac{1}{2}, s) = 2\gamma_1(2x, s) - \gamma_1(x, s)$$

$$+ 2^{1-s} \log^2 2 \cdot C_s(2x) + 2^{2-s} \log 2 \sum_{n=1}^{\infty} \frac{\log n \cos(4n\pi x)}{n^s}$$

$$- 2^{1-s} \pi \log 2 \cdot S_s(2x) + 2^{2-s} \log 2 [\gamma + \log(2\pi)] C_s(2x)$$

We then compute the limit $s \to 0$

$$\gamma_1(x + \tfrac{1}{2}) = 2\gamma_1(2x) - \gamma_1(x) - \log^2 2 + 4 \log 2 \lim_{s \to 0} \sum_{n=1}^{\infty} \frac{\log n \cos(4n\pi x)}{n^s}$$

$$- \pi \log 2 \cdot \cot 2\pi x - 2\log 2[\gamma + \log(2\pi)]$$

since $\lim_{s \to 0} C_s(x) = -\frac{1}{2}$ and $\lim_{s \to 0} S_s(x) = \frac{1}{2} \cot \pi x$.

Using (2.8) we obtain the desired result

(7.20)  $\gamma_1(x + \tfrac{1}{2}) = 2\gamma_1(2x) - \gamma_1(x) - \log^2 2 + 2\log 2 \cdot \psi(2x)$

For example, with $x = \tfrac{1}{4}$ in (7.20) we have

$$\gamma_1\left(\frac{1}{4}\right)+\gamma_1\left(\frac{3}{4}\right)=2\gamma_1\left(\frac{1}{2}\right)-\log^2 2+2\log 2\cdot\psi\left(\frac{1}{2}\right)$$

and employing

$$\gamma_1\left(\frac{1}{2}\right)=\gamma_1-\log^2 2-2\gamma\log 2$$

$$\psi\left(\frac{1}{2}\right)=-\gamma-2\log 2$$

we obtain

$$\gamma_1\left(\frac{1}{4}\right)+\gamma_1\left(\frac{3}{4}\right)=2\gamma_1-7\log^2 2-6\gamma\log 2$$

as may be seen above.

Another derivation of (7.20) follows below.

$\square$

Hansen and Patrick [33] showed in 1962 that the Hurwitz zeta function could be written as

(7.21) $\qquad \varsigma\left(s,x+\frac{1}{2}\right)=2^s\varsigma(s,2x)-\varsigma(s,x)$

and, by analytic continuation, this holds for all $s$. This may also be obtained by letting $x\to x+\frac{1}{2}$ in (2.1) which gives us

(7.22) $\qquad \varsigma(s,x+\tfrac{1}{2})=2\Gamma(1-s)\left[\sin\left(\frac{\pi s}{2}\right)\sum_{n=1}^{\infty}(-1)^n\frac{\cos 2n\pi x}{(2\pi n)^{1-s}}+\cos\left(\frac{\pi s}{2}\right)\sum_{n=1}^{\infty}(-1)^n\frac{\sin 2n\pi x}{(2\pi n)^{1-s}}\right]$

where $\mathrm{Re}(s)<0$ and $-\frac{1}{2}<x\le\frac{1}{2}$. To deduce (7.21) we then apply the following relationship for suitably convergent series

$$\sum_{n=1}^{\infty}(-1)^n a_n=-\sum_{n=1}^{\infty}a_n+2\sum_{n=1}^{\infty}a_{2n}$$

We then multiply (7.21) by $s-1$

$$(s-1)\varsigma\left(s,x+\frac{1}{2}\right)=2^s(s-1)\varsigma(s,2x)-(s-1)\varsigma(s,x)$$

Differentiation gives us

$$\frac{\partial}{\partial s}\left[(s-1)\varsigma\left(s,x+\frac{1}{2}\right)\right]=2^s\log 2.(s-1)\varsigma(s,2x)+2^s\frac{\partial}{\partial s}[(s-1)\varsigma(s,2x)]-\frac{\partial}{\partial s}[(s-1)\varsigma(s,x)]$$

and, with $s=1$, we immediately obtain (3.6) using (1.2.2).

The second derivative gives us

$$\frac{\partial^2}{\partial s^2}\left[(s-1)\varsigma\left(s,x+\frac{1}{2}\right)\right]=2^{s+1}\log 2.\frac{\partial}{\partial s}(s-1)\varsigma(s,2x)+2^s\log^2 2.(s-1)\varsigma(s,2x)$$

$$+2^s\frac{\partial^2}{\partial s^2}[(s-1)\varsigma(s,2x)]-\frac{\partial^2}{\partial s^2}[(s-1)\varsigma(s,x)]$$

so that with $s=1$ we immediately obtain another derivation of the functional equation (7.20).

A third derivation of (7.20) is presented in Section 8 below.

$\square$

It is interesting to note that substituting Hurwitz's formula (1.3) in (7.10) gives us

$$\varsigma_a(s,t)=2^{1-s}\Gamma(1-s)\left[\sin\left(\frac{\pi s}{2}\right)\sum_{n=1}^{\infty}\frac{\cos n\pi t-\cos n\pi(1+t)}{(2\pi n)^{1-s}}+\cos\left(\frac{\pi s}{2}\right)\sum_{n=1}^{\infty}\frac{\sin n\pi t-\sin n\pi(1+t)}{(2\pi n)^{1-s}}\right]$$

$$=2^{1-s}\Gamma(1-s)\left[\sin\left(\frac{\pi s}{2}\right)\sum_{n=1}^{\infty}\frac{[1-(-1)^n]\cos n\pi t}{(2\pi n)^{1-s}}+\cos\left(\frac{\pi s}{2}\right)\sum_{n=1}^{\infty}\frac{[1-(-1)^n]\sin n\pi t}{(2\pi n)^{1-s}}\right]$$

and we obtain

(7.23) $$\varsigma_a(s,t)=2\Gamma(1-s)\pi^{1-s}\left[\sin\left(\frac{\pi s}{2}\right)\sum_{n=0}^{\infty}\frac{\cos(2n+1)\pi t}{(2n+1)^{1-s}}+\cos\left(\frac{\pi s}{2}\right)\sum_{n=0}^{\infty}\frac{\sin(2n+1)\pi t}{(2n+1)^{1-s}}\right]$$

By means of contour integration, Boudjelkha [12] showed that (7.23) was valid for $\sigma<0, 0<t\leq 1;\ 0<\sigma, t<1$.

**8. The Fourier Series for $\varsigma''(0,x)$**

**Theorem**

We have for $0<u<1$

(8.1) $$\varsigma''(0,u)=\sum_{n=1}^{\infty}\frac{[\gamma+\log(2n\pi)]^2\sin(2n\pi u)}{n\pi}-\frac{1}{4}\varsigma(2)\sum_{n=1}^{\infty}\frac{\sin(2n\pi u)}{n\pi}$$

$$+\sum_{n=1}^{\infty}\frac{[\gamma+\log(2\pi n)]\cos(2n\pi u)}{n}$$

**Proof**

Assuming that for $0 < u < 1$

$$\lim_{s \to 0} \int_{\frac{1}{2}}^{u} \gamma_1(x, s) \, dx = \int_{\frac{1}{2}}^{u} \lim_{s \to 0} \gamma_1(x, s) \, dx$$

$$= \int_{\frac{1}{2}}^{u} \gamma_1(x) \, dx$$

we have

$$\lim_{s \to 0} \int_{\frac{1}{2}}^{u} \gamma_1(x, s) \, dx = \sum_{n=1}^{\infty} \frac{\log^2 n \sin(2n\pi u)}{2n\pi} + 2[\gamma + \log(2\pi)] \sum_{n=1}^{\infty} \frac{\log n \sin(2n\pi u)}{2n\pi}$$

$$+ \left[ [\gamma + \log(2\pi)]^2 - \frac{1}{2}\varsigma(2) \right] \sum_{n=1}^{\infty} \frac{\sin(2n\pi u)}{2n\pi}$$

$$- \pi \sum_{n=1}^{\infty} \frac{\log n \left[ (-1)^n - \cos(2n\pi u) \right]}{2n\pi} - \pi [\gamma + \log(2\pi)] \sum_{n=1}^{\infty} \frac{\left[ (-1)^n - \cos(2n\pi u) \right]}{2n\pi}$$

As shown by Dilcher [30] we have for $k \geq 0$

$$(-1)^{k+1} \varsigma_a^{(k)}(1) = \sum_{n=1}^{\infty} (-1)^n \frac{\log^k n}{n}$$

$$= \sum_{j=0}^{k-1} \binom{k}{j} \gamma_j \log^{k-j} 2 - \frac{\log^{k+1} 2}{k+1}$$

and, in particular, we have

$$\sum_{n=1}^{\infty} (-1)^n \frac{\log n}{n} = \gamma \log 2 - \frac{1}{2} \log^2 2$$

This gives us

$$\lim_{s \to 0} \int_{\frac{1}{2}}^{u} \gamma_1(x, s) \, dx = \sum_{n=1}^{\infty} \frac{\log^2 n \sin(2n\pi u)}{2n\pi} + [\gamma + \log(2\pi)] \sum_{n=1}^{\infty} \frac{\log n \sin(2n\pi u)}{n\pi}$$

$$+ \left[ [\gamma + \log(2\pi)]^2 - \frac{1}{2}\varsigma(2) \right] \sum_{n=1}^{\infty} \frac{\sin(2n\pi u)}{2n\pi}$$

$$- \frac{1}{2} \left[ \gamma \log 2 - \frac{1}{2} \log^2 2 \right] + \frac{1}{2} \sum_{n=1}^{\infty} \frac{[\gamma + \log(2\pi n)] \cos(2n\pi u)}{n}$$

$$+\frac{1}{2}[\gamma+\log(2\pi)]\log 2$$

It is known that [27]

(8.2) $$\gamma_n(x) = \frac{(-1)^{n+1}}{n+1} \frac{\partial}{\partial x} \varsigma^{(n+1)}(0,x)$$

and thus

(8.3) $$\int_1^u \gamma_n(x)\, dx = \frac{(-1)^{n+1}}{n+1}\left[\varsigma^{(n+1)}(0,u) - \varsigma^{(n+1)}(0)\right]$$

Hence we have

$$\int_1^u \gamma_1(x)\, dx = \frac{1}{2}\left[\varsigma''(0,u) - \varsigma''(0)\right]$$

Since

$$\int_{\frac{1}{2}}^u \gamma_1(x)\, dx = \int_1^u \gamma_1(x)\, dx - \int_1^{\frac{1}{2}} \gamma_1(x)\, dx$$

we obtain

$$\int_{\frac{1}{2}}^u \gamma_1(x)\, dx = \frac{1}{2}\left[\varsigma''(0,u) - \varsigma''\!\left(0,\tfrac{1}{2}\right)\right]$$

Using the well-known formula $\varsigma(s,\tfrac{1}{2}) = (2^s - 1)\varsigma(s)$ we readily find that

$$\varsigma''\!\left(0,\tfrac{1}{2}\right) = 2\varsigma'(0)\log 2 - \frac{1}{2}\log^2 2$$

and Lerch's formula

$$\varsigma'(0,x) = \log\Gamma(x) - \frac{1}{2}\log(2\pi)$$

gives us

$$\varsigma''\!\left(0,\tfrac{1}{2}\right) = -\log(2\pi)\log 2 - \frac{1}{2}\log^2 2$$

Accordingly, with a little algebra we obtain

$$\varsigma''(0,u) = \sum_{n=1}^{\infty} \frac{[\gamma+\log(2n\pi)]^2 \sin(2n\pi u)}{n\pi} - \frac{1}{4}\varsigma(2)\sum_{n=1}^{\infty} \frac{\sin(2n\pi u)}{n\pi}$$

$$+ \sum_{n=1}^{\infty} \frac{[\gamma+\log(2\pi n)]\cos(2n\pi u)}{n}$$

This was in fact obtained much more directly in [24] by differentiating (2.1) twice with respect to $s$ and setting $s = 0$.

We easily determine from the above that

$$\int_0^1 \varsigma''(0,u)\,du = 0$$

as previously reported by Deninger [28].

Letting $u \to 1-u$ we find that

(8.4) $$\sum_{n=1}^{\infty} \frac{\gamma + \log(2\pi n)}{n} \cos(2n\pi u) = \frac{1}{2}[\varsigma''(0,u) + \varsigma''(0,1-u)]$$

This identity was derived by Ramanujan and it appears, albeit in a highly disguised form, in Berndt's book [4, Part I, p.208] and [5]. This identity was reported by Deninger [28] in 1984 and was also derived by Blagouchine [11] in 2014.

We also have

(8.5) $$\frac{1}{2}[\varsigma''(0,u) - \varsigma''(0,1-u)] = \sum_{n=1}^{\infty} \frac{[\gamma + \log(2n\pi)]^2 \sin(2n\pi u)}{n\pi} - \frac{1}{4}\varsigma(2)\sum_{n=1}^{\infty} \frac{\sin(2n\pi u)}{n\pi}$$

Using (4.2) we may write (8.4) as

(8.6) $$\frac{1}{2}[\varsigma''(0,u) + \varsigma''(0,1-u)] = \sum_{n=1}^{\infty} \frac{\log n}{n} \cos(2n\pi u) - [\gamma + \log(2\pi)]\log[2\sin \pi u]$$

With $u = \tfrac{1}{3}$ in (8.6) we obtain

$$\frac{1}{2}[\varsigma''(0,\tfrac{1}{3}) + \varsigma''(0,\tfrac{2}{3})] = \sum_{n=1}^{\infty} \frac{\log n}{n} \cos(2n\pi/3) - \frac{1}{2}[\gamma + \log(2\pi)]\log 3$$

which was previously derived by Coffey [17].

It is well known [33] that

(8.7) $$\sum_{r=1}^{q-1} \varsigma\left(s, \frac{r}{q}\right) = (q^s - 1)\varsigma(s)$$

and hence

$$\sum_{r=1}^{q-1} \varsigma''\left(s, \frac{r}{q}\right) = (q^s - 1)\varsigma''(s) + 2q^s \varsigma'(s)\log q + q^s \varsigma(s)\log^2 q$$

This gives us

$$\varsigma''\left(0,\frac{1}{3}\right)+\varsigma''\left(0,\frac{2}{3}\right)=2\varsigma'(0)\log 3+\varsigma(0)\log^2 3$$

$$=-\log(2\pi)\log 3-\frac{1}{2}\log^2 3$$

Hence we obtain

(8.8) $$\sum_{n=1}^{\infty}\frac{\log n}{n}\cos(2n\pi/3)=\frac{1}{2}\left[\gamma-\frac{1}{2}\log 3\right]\log 3$$

as previously derived by Candelpergher [13, p.155] in a different manner.

It was mentioned above in (1.7) that Srivastava and Tsumura [45, p.293] reported for $\mathrm{Re}(s)>1$

$$\sum_{n=1}^{\infty}\frac{\cos(2n\pi/3)}{n^s}=\frac{1}{2}(3^{1-s}-1)\varsigma(s)$$

which may be expressed as

(8.9) $$\sum_{n=1}^{\infty}\frac{\cos(2n\pi/3)}{n^s}=\frac{1}{2}\frac{3^{1-s}-1}{s-1}(s-1)\varsigma(s)$$

Since $x^{1-s}-1=\exp[-(s-1)\log x]-1$ we have

$$\frac{3^{1-s}-1}{s-1}=\sum_{n=1}^{\infty}\frac{(-1)^n\log^n 3}{n!}(s-1)^{n-1}$$

Differentiation of (8.9) results in

$$-2\sum_{n=1}^{\infty}\frac{\log n}{n^s}\cos(2n\pi/3)=(s-1)\varsigma(s)\sum_{n=1}^{\infty}\frac{(-1)^n(n-1)\log^n 3}{n!}(s-1)^{n-2}$$

$$+\sum_{n=1}^{\infty}\frac{(-1)^n\log^n 3}{n!}(s-1)^{n-1}\frac{d}{ds}[(s-1)\varsigma(s)]$$

and with $s=1$ we thereby obtain (8.8) by a different method.

Another differentiation gives us

$$2\sum_{n=1}^{\infty}\frac{\log^2 n}{n^s}\cos(2n\pi/3)=(s-1)\varsigma(s)\sum_{n=1}^{\infty}\frac{(-1)^n(n-1)(n-2)\log^n 3}{n!}(s-1)^{n-3}$$

$$+2\sum_{n=1}^{\infty}\frac{(-1)^n(n-1)\log^n 3}{n!}(s-1)^{n-2}\frac{d}{ds}[(s-1)\varsigma(s)]$$

$$+\sum_{n=1}^{\infty}\frac{(-1)^n\log^n 3}{n!}(s-1)^{n-1}\frac{d^2}{ds^2}[(s-1)\varsigma(s)]$$

and with $s=1$ we obtain

(8.10) $\quad 2\sum_{n=1}^{\infty}\frac{\log^2 n}{n}\cos(2n\pi/3)=-\frac{1}{3}\log^3 3+\gamma\log^2 3+2\gamma_1\log 2$

which was also derived by Candelpergher [13, p.155]. It should be noted that it has not been proved that (1.7) actually holds in the case $s=1$.

**Theorem**

We have for $0<x<\frac{1}{2}$

(8.11) $\quad \varsigma''(0,x+\tfrac{1}{2})=\varsigma''(0,2x)-\varsigma''(0,x)+2\log 2\left[\log\Gamma(2x)-\frac{1}{2}\log(2\pi)\right]+\log^2 2\left[\frac{1}{2}-2x\right]$

**Proof**

We let $x\to x+\tfrac{1}{2}$ in (8.1) to obtain for $|x|<\tfrac{1}{2}$

$$\varsigma''(0,x+\tfrac{1}{2})=\sum_{n=1}^{\infty}(-1)^n\frac{[\gamma+\log(2n\pi)]^2\sin(2n\pi x)}{n\pi}-\frac{1}{4}\varsigma(2)\sum_{n=1}^{\infty}(-1)^n\frac{\sin(2n\pi x)}{n\pi}$$

$$+\sum_{n=1}^{\infty}(-1)^n\frac{[\gamma+\log(2\pi n)]\cos(2n\pi x)}{n}$$

$$=\sum_{n=1}^{\infty}\frac{[\gamma+\log(4n\pi)]^2\sin(4n\pi x)}{n\pi}-\frac{1}{4}\varsigma(2)\sum_{n=1}^{\infty}\frac{\sin(4n\pi x)}{n\pi}$$

$$+\sum_{n=1}^{\infty}\frac{[\gamma+\log(4\pi n)]\cos(4n\pi x)}{n}-\varsigma''(0,x)$$

where we have employed (3.1.1).

$$=\sum_{n=1}^{\infty}\frac{[\gamma+\log(2n\pi)]^2\sin(4n\pi x)}{n\pi}+2\log 2\sum_{n=1}^{\infty}\frac{[\gamma+\log(2n\pi)]\sin(4n\pi x)}{n\pi}$$

$$+\left[\log^2 2-\frac{1}{4}\varsigma(2)\right]\sum_{n=1}^{\infty}\frac{\sin(4n\pi x)}{n\pi}$$

$$+\sum_{n=1}^{\infty}\frac{[\gamma+\log(2\pi n)]\cos(4n\pi x)}{n}+\log 2\sum_{n=1}^{\infty}\frac{\cos(4n\pi x)}{n}-\varsigma''(0,x)$$

$$= \varsigma''(0, 2x) + 2\log 2 \sum_{n=1}^{\infty} \frac{[\gamma + \log(2n\pi)]\sin(4n\pi x)}{n\pi}$$

$$+ \log^2 2 \sum_{n=1}^{\infty} \frac{\sin(4n\pi x)}{n\pi} + \log 2 \sum_{n=1}^{\infty} \frac{\cos(4n\pi x)}{n} - \varsigma''(0, x)$$

We then employ (5.2), (5.3) and (5.5) for $0 < x < \frac{1}{2}$ to obtain

(8.12) $\quad \varsigma''\left(0, x + \tfrac{1}{2}\right) = \varsigma''(0, 2x) - \varsigma''(0, x) + 2\log 2 \left[\log \Gamma(2x) - \frac{1}{2}\log(2\pi)\right] + \log^2 2 \left[\frac{1}{2} - 2x\right]$

as was also shown in [24] by a different method. It is easily seen that this is in fact valid for $x = \tfrac{1}{2}$ where we see that

$$\varsigma''\left(0, \tfrac{1}{2}\right) = -\log(2\pi)\log 2 - \frac{1}{2}\log^2 2$$

As shown below, the functional equation (8.10) is actually valid for all $x > 0$. To see this, we differentiate

$$\varsigma\left(s, x + \frac{1}{2}\right) = 2^s \varsigma(s, 2x) - \varsigma(s, x)$$

so that

$$\varsigma'\left(s, x + \tfrac{1}{2}\right) = 2^s \log 2 . \varsigma(s, 2x) + 2^s \varsigma'(s, 2x) - \varsigma'(s, x)$$

The second derivative gives us

$$\varsigma''\left(s, x + \tfrac{1}{2}\right) = 2^s \log^2 2 . \varsigma(s, 2x) + 2^{s+1} \log 2 . \varsigma'(s, 2x) + 2^s \varsigma''(s, 2x) - \varsigma''(s, x)$$

so that with $s = 0$ we immediately obtain another derivation of the functional equation (8.10)

(8.13) $\quad \varsigma''\left(0, x + \tfrac{1}{2}\right) = \log^2 2 . \varsigma(0, 2x) + 2\log 2 . \varsigma'(0, 2x) + \varsigma''(0, 2x) - \varsigma''(0, x)$

or

(8.14) $\quad \varsigma''\left(0, x + \tfrac{1}{2}\right) = \left(\frac{1}{2} - 2x\right) \log^2 2 + 2\log 2 \left[\log \Gamma(2x) - \frac{1}{2}\log(2\pi)\right] + \varsigma''(0, 2x) - \varsigma''(0, x)$

$\square$

We know from (8.2) that

$$\gamma_1(x) = \frac{1}{2} \frac{\partial}{\partial x} \varsigma''(0, x)$$

and differentiating (8.11) therefore results in

$$\gamma_1(x + \tfrac{1}{2}) = 2\gamma_1(2x) - \gamma_1(x) - \log^2 2 + 2\log 2 \cdot \psi(2x)$$

which we saw above in (7.19).

□

Reference to (1.2.2) shows that

$$f(s)\varsigma(s,x) - g(s)\varsigma(s,y) = \frac{f(s) - g(s)}{s-1} + \sum_{n=0}^{\infty} \frac{(-1)^n}{n!}[f(s)\gamma_n(x) - g(s)\gamma_n(y)](s-1)^n$$

and, assuming that $f(1) = g(1)$, L'Hôpital's rule gives us

$$\lim_{s \to 1}[f(s)\varsigma(s,x) - g(s)\varsigma(s,y)] = f'(1) - g'(1) - f(1)[\psi(x) - \psi(y)]$$

Differentiation results in

$$\frac{d}{ds}[f(s)\varsigma(s,x) - g(s)\varsigma(s,y)] = \frac{d}{ds}\frac{f(s) - g(s)}{s-1}$$

$$+ \sum_{n=0}^{\infty} \frac{(-1)^n}{n!}[f(s)\gamma_n(x) - g(s)\gamma_n(y)]n(s-1)^{n-1}$$

$$+ \sum_{n=0}^{\infty} \frac{(-1)^n}{n!}[f'(s)\gamma_n(x) - g'(s)\gamma_n(y)](s-1)^n$$

and accordingly we have

$$\lim_{s \to 1}\frac{d}{ds}[f(s)\varsigma(s,x) - g(s)\varsigma(s,y)] = \frac{d}{ds}\frac{f(s) - g(s)}{s-1}\bigg|_{s=1} - f(1)[\gamma_1(x) - \gamma_1(y)]$$

$$-[f'(1)\psi(x) - g'(1)\psi(y)]$$

We write (1.11) as

$$\sum_{n=1}^{\infty} \frac{\sin(n\pi/2)}{n^s} = 2^{1-2s}\varsigma\left(s, \frac{1}{4}\right) - (1 - 2^{-s})\varsigma(s)$$

$$:= f(s)\varsigma\left(s, \frac{1}{4}\right) - g(s)\varsigma(s)$$

where

$$f(s) = 2^{1-2s} \qquad\qquad g(s) = 1 - 2^{-s}$$

We see that

$$f(1) = \frac{1}{2} \qquad\qquad g(1) = \frac{1}{2}$$

$$f'(s) = -2^{2(1-s)} \log 2 \qquad\qquad g'(s) = 2^{-s} \log 2$$

$$f'(1) = -\log 2 \qquad\qquad g'(1) = \frac{1}{2}\log 2$$

Simple algebra gives us

$$f(s) - g(s) = 2^{1-2s} + 2^{-s} - 1$$

$$= \frac{1}{2} 4^{1-s} + \frac{1}{2} 2^{1-s} - 1$$

$$= \frac{1}{2}[4^{1-s} + 2^{1-s} - 2]$$

$$= \frac{1}{2}[e^{(1-s)\log 4} + e^{(1-s)\log 2} - 2]$$

$$= \frac{1}{2} \sum_{n=1}^{\infty} \frac{(-1)^n}{n!}[\log^n 4 + \log^n 2](s-1)^n$$

Therefore we have

$$\frac{f(s) - g(s)}{s-1} = \frac{1}{2}\sum_{n=1}^{\infty} \frac{(-1)^n}{n!}[\log^n 4 + \log^n 2](s-1)^{n-1}$$

and thus

$$\frac{d}{ds} \frac{f(s)-g(s)}{s-1}\bigg|_{s=1} = \frac{1}{4}[\log^2 2 + \log^2 4]$$

Hence we obtain

$$\sum_{n=1}^{\infty} \frac{\log n}{n} \sin(n\pi/2) = -\frac{1}{4}[\log^2 2 + \log^2 4] + f(1)[\gamma_1(\tfrac{1}{4}) - \gamma_1] + f'(1)\psi(\tfrac{1}{4}) - g'(1)\psi(1)$$

$$= -\frac{5}{4}\log^2 2 + \frac{1}{2}[\gamma_1(\tfrac{1}{4}) - \gamma_1] - \psi(\tfrac{1}{4})\log 2 + \frac{1}{2}\gamma \log 2$$

Employing $\psi(\tfrac{1}{4}) = -\gamma - \frac{1}{2}\pi - 3\log 2$ and Kummer's formula (5.4)

$$\frac{1}{\pi}\sum_{n=1}^{\infty} \frac{\log n}{n} \sin(n\pi/2) = \log \Gamma(\tfrac{1}{4}) - \frac{1}{4}\gamma - \frac{3}{4}\log \pi - \frac{1}{2}\log 2$$

we obtain another derivation of the known result

$$\gamma_1\left(\frac{1}{4}\right) = \frac{1}{2}[2\gamma_1 - 7\log^2 2 - 6\gamma \log 2] - \frac{1}{2}\pi\left[\gamma + 4\log 2 + 3\log \pi - 4\log\Gamma\left(\frac{1}{4}\right)\right]$$

Taking the second derivative, it is apparent that $\sum_{n=1}^{\infty}\frac{\log^2 n}{n}\sin(n\pi/2)$ will involve $\gamma_2(\tfrac{1}{4}) - \gamma_2$, whereupon (8.1) will be germane.

## 9. Analysis provided by Bernard Candelpergher

Shortly after I submitted my earlier paper [27] on this subject, I contacted Bernard Candelpergher to briefly point out the very close connection between some of the formulae appearing in that paper with those contained in his work [14] regarding Ramanujan summation. He kindly provided me with the following information by email in October 2019:

Candelpergher's book [14, p.5] reports that if the series $\sum_{n\geq 1} f(n)$ and the integral $\int_1^{\infty} f(x)\,dx$ are convergent then we have

$$(9.1) \qquad \sum_{n\geq 1}^{R} f(n) = \sum_{n\geq 1}^{\infty} f(n) - \int_1^{\infty} f(t)\,dt$$

where $\sum_{n\geq 1}^{R} f(n)$ is the Ramanujan constant.

Since $\sum_{n=1}^{\infty}\frac{e^{2\pi inx}}{n^s}$ is convergent for $\operatorname{Re} s > 1$ we have

$$\sum_{n\geq 1}^{R}\frac{e^{2\pi inx}}{n^s} = \sum_{n\geq 1}^{\infty}\frac{e^{2\pi inx}}{n^s} - \int_1^{\infty}\frac{e^{2\pi itx}}{t^s}\,dt$$

Integration by parts gives us

$$\int_1^{\infty}\frac{e^{2\pi itx}}{t^s}\,dt = \left.\frac{e^{2\pi itx}}{2\pi ixt^s}\right|_1^{\infty} + \frac{s}{2\pi ix}\int_1^{\infty}\frac{e^{2\pi itx}}{t^{s+1}}\,dt$$

and we have the limit

$$\lim_{s\to 0}\int_1^{\infty}\frac{e^{2\pi itx}}{t^s}\,dt = -\frac{e^{2\pi ix}}{2\pi ix}$$

Therefore

(9.2) $$\lim_{s \to 0} \sum_{n \geq 1}^{R} \frac{e^{2\pi i n x}}{n^s} = \lim_{s \to 0} \sum_{n \geq 1}^{\infty} \frac{e^{2\pi i n x}}{n^s} + \frac{e^{2\pi i x}}{2\pi i x}$$

We have

$$\lim_{s \to 0} \sum_{n \geq 1}^{R} \frac{e^{2\pi i n x}}{n^s} = \sum_{n \geq 1}^{R} e^{2\pi i n x}$$

and, applying (9.1) again, we have

$$\sum_{n \geq 1}^{R} e^{2\pi i n x} = \sum_{n=1}^{\infty} e^{2\pi i n x} - \int_{1}^{\infty} e^{2\pi i t x} dt$$

Therefore we have

(9.3) $$\lim_{s \to 0} \sum_{n \geq 1}^{R} \frac{e^{2\pi i n x}}{n^s} = \frac{e^{2\pi i n x}}{1 - e^{2\pi i n x}} + \frac{e^{2\pi i x}}{2\pi i x}$$

Equating (9.2) and (9.3) results in

(9.4) $$\lim_{s \to 0} \sum_{n \geq 1}^{\infty} \frac{e^{2\pi i n x}}{n^s} = \frac{e^{2\pi i n x}}{1 - e^{2\pi i n x}}$$

in accordance with (1.4) and (1.5).

$\square$

Candelpergher also considered the limit $\lim_{s \to 0} \sum_{n \geq 1}^{\infty} \log n \frac{e^{2\pi i n x}}{n^s}$ and, since this is convergent for $\operatorname{Re} s > 1$, we have in accordance with (9.1)

$$\sum_{n \geq 1}^{R} \log n \frac{e^{2\pi i n x}}{n^s} = \sum_{n \geq 1}^{\infty} \log n \frac{e^{2\pi i n x}}{n^s} - \int_{1}^{\infty} \log t \frac{e^{2\pi i t x}}{t^s} dt$$

Integration by parts gives us

$$\int_{1}^{\infty} \frac{\log t \, e^{2\pi i t x}}{t^s} dt = \frac{\log t \, e^{2\pi i t x}}{2\pi i x t^s} \bigg|_{1}^{\infty} - \frac{1}{2\pi i x} \int_{1}^{\infty} \frac{(1 - s \log t) e^{2\pi i t x}}{t^{s+1}} dt$$

and we have the limit

$$\lim_{s \to 0} \int_{1}^{\infty} \log t \frac{e^{2\pi i t x}}{t^s} dt = -\frac{1}{2\pi i x} \int_{1}^{\infty} \frac{e^{2\pi i t x}}{t} dt$$

Hence we have

(9.5) $$\lim_{s \to 0} \sum_{n \geq 1}^{R} \log n \frac{e^{2\pi i n x}}{n^s} = \lim_{s \to 0} \sum_{n \geq 1}^{\infty} \log n \frac{e^{2\pi i n x}}{n^s} + \frac{1}{2\pi i x} \int_{1}^{\infty} \frac{e^{2\pi i t x}}{t} dt$$

We see that

$$\lim_{s \to 0} \sum_{n \geq 1}^{R} \log n \frac{e^{2\pi i n x}}{n^s} = \sum_{n \geq 1}^{R} \log n \cdot e^{2\pi i n x}$$

Candelpergher et al. ([13] and [14, p.68]) have shown that

(9.6) $$\sum_{n \geq 1}^{R} \log n \cdot e^{2\pi i n x} = \sum_{n=0}^{\infty} \frac{(-1)^{n-1}}{n!} \varsigma'(-n)(-2\pi i x)^n - \frac{1}{2\pi i x} \int_{0}^{1} \frac{e^{2\pi i t x} - 1}{t} dt$$

The real part of (9.6) is

(9.7) $$\sum_{n \geq 1}^{R} \log n \cdot \cos 2\pi n x = \sum_{n=0}^{\infty} \frac{(-1)^{n-1}}{(2n)!} \varsigma'(-2n)(2\pi x)^{2n} - \frac{1}{2\pi x} \int_{0}^{1} \frac{\sin 2\pi x t}{t} dt$$

Noting that

$$\varsigma'(0) = -\frac{1}{2} \log(2\pi)$$

and [44, p.99]

$$\varsigma'(-2n) = (-1)^n \frac{(2n)!}{2(2\pi)^{2n}} \varsigma(2n+1)$$

and

$$\frac{1}{2\pi x} \int_{0}^{1} \frac{\sin 2\pi x t}{t} dt = \frac{1}{2\pi x} \int_{0}^{\infty} \frac{\sin 2\pi x t}{t} dt - \frac{1}{2\pi x} \int_{1}^{\infty} \frac{\sin 2\pi x t}{t} dt$$

$$= \frac{1}{4x} - \frac{1}{2\pi x} \int_{1}^{\infty} \frac{\sin 2\pi x t}{t} dt$$

and using the well-known integral from Fourier series analysis [3, p.286]

$$\frac{\pi}{2} = \int_{0}^{\infty} \frac{\sin t}{t} dt$$

we may write (9.7) as

$$\sum_{n \geq 1}^{R} \log n \cdot \cos 2\pi n x = \frac{1}{2} \log(2\pi) - \frac{1}{2} \sum_{n=1}^{\infty} \varsigma(2n+1) x^{2n} - \frac{1}{4x} + \frac{1}{2\pi x} \int_{1}^{\infty} \frac{\sin 2\pi x t}{t} dt$$

It is known [44, p.160] that for $|x| < 1$

$$\sum_{n=1}^{\infty} \varsigma(2n+1)x^{2n} = -\frac{1}{2}[\psi(1+x) + \psi(1-x)] - \gamma$$

and hence we obtain

(9.8) $\quad \sum_{n\geq 1}^{R} \log n . \cos 2\pi nx = \frac{1}{4}[\psi(1+x) + \psi(1-x)] + \frac{1}{2}[\gamma + \log(2\pi)] - \frac{1}{4x} + \frac{1}{2\pi x}\int_{1}^{\infty} \frac{\sin 2\pi xt}{t} dt$

or equivalently

$$= \frac{1}{4}[\psi(x) + \psi(1-x)] + \frac{1}{2}[\gamma + \log(2\pi)] + \frac{1}{2\pi x}\int_{1}^{\infty} \frac{\sin 2\pi xt}{t} dt$$

$$= \frac{1}{2}\left[\psi(x) + \frac{\pi}{2}\cot \pi x + \gamma + \log(2\pi)\right] + \frac{1}{2\pi x}\int_{1}^{\infty} \frac{\sin 2\pi xt}{t} dt$$

Comparing this with (9.5) we obtain

(9.9) $\quad \lim_{s\to 0}\sum_{n\geq 1}^{\infty} \frac{\log n}{n^s}\cos 2\pi nx = \frac{1}{4}[\psi(x) + \psi(1-x)] + \frac{1}{2}[\gamma + \log(2\pi)]$

which concurs with (2.8).

$\square$

The imaginary part of (9.6) is

(9.10) $\quad \sum_{n\geq 1}^{R} \log n . \sin 2\pi nx = \sum_{n=0}^{\infty} \frac{(-1)^{n-1}}{(2n+1)!}\varsigma'(-2n-1)(2\pi x)^{2n+1} + \frac{1}{2\pi x}\int_{0}^{1} \frac{\cos 2\pi xt - 1}{t} dt$

and from (9.5) we have

(9.11) $\quad \sum_{n\geq 1}^{R} \log n . \sin 2\pi nx = \lim_{s\to 0}\sum_{n\geq 1}^{\infty} \log n \frac{\sin 2\pi nx}{n^s} - \frac{1}{2\pi x}\int_{1}^{\infty} \frac{\cos 2\pi xt}{t} dt$

Therefore we have

$$\lim_{s\to 0}\sum_{n\geq 1}^{\infty} \log n \frac{\sin 2\pi nx}{n^s} = \sum_{n=0}^{\infty} \frac{(-1)^{n-1}}{(2n+1)!}\varsigma'(-2n-1)(2\pi x)^{2n+1} + \frac{1}{2\pi x}\int_{0}^{1} \frac{\cos 2\pi xt - 1}{t} dt$$

$$+ \frac{1}{2\pi x}\int_{1}^{\infty} \frac{\cos 2\pi nx}{t} dt$$

Using

$$\int\limits_0^1 \frac{\cos 2\pi xt - 1}{t}\,dt = \int\limits_0^{2\pi x} \frac{\cos t - 1}{t}\,dt$$

$$\int\limits_1^\infty \frac{\cos 2\pi xt}{t}\,dt = \int\limits_{2\pi x}^\infty \frac{\cos t}{t}\,dt$$

and the known formula [23] for the cosine integral function $Ci(x)$

$$Ci(x) = -\int\limits_x^\infty \frac{\cos t}{t}\,dt = \gamma + \log x + \int\limits_0^x \frac{\cos t - 1}{t}\,dt$$

we see that

$$\int\limits_0^{2\pi x} \frac{\cos t - 1}{t}\,dt + \int\limits_{2\pi x}^\infty \frac{\cos t}{t}\,dt = -[\gamma + \log(2\pi x)]$$

This results in

(9.12) $$\lim_{s \to 0} \sum_{n \geq 1} \log n \frac{\sin 2\pi nx}{n^s} = \sum_{n=0}^\infty \frac{(-1)^{n-1}}{(2n+1)!} \varsigma'(-2n-1)(2\pi x)^{2n+1} - \frac{\gamma + \log(2\pi x)}{2\pi x}$$

We saw in (1.1) that

$$\gamma_1(1-x, s) - \gamma_1(x, s) = 2\pi \sum_{n=1}^\infty \frac{[\gamma + \log(2\pi n)]}{n^s} \sin(2n\pi x)$$

and in the limit as $s \to 0$ we have

(9.13) $$\gamma_1(1-x) - \gamma_1(x) = 2\pi \lim_{s \to 1} \sum_{n=1}^\infty \frac{\sin(2n\pi x) \log n}{n^{1-s}} + \pi[\gamma + \log(2\pi)]\cot(\pi x)$$

which concurs with Proposition 4.2 in [27].

Accordingly, assuming the veracity of (9.13), we deduce that

(9.14) $$\frac{1}{2\pi}[\gamma_1(1-x) - \gamma_1(x)] = \sum_{n=0}^\infty \frac{(-1)^{n-1}}{(2n+1)!} \varsigma'(-2n-1)(2\pi x)^{2n+1} - \frac{\gamma + \log(2\pi x)}{2\pi x}$$

$$+ \frac{1}{2}[\gamma + \log(2\pi)]\cot(\pi x)$$

**Comments by the author:**

Using

$$\gamma_n(1+x) - \gamma_n(x) = -\frac{\log^n x}{x}$$

and

$$\psi(1+x) - \psi(1-x) = -\pi \cot(\pi x) - \frac{1}{x}$$

we may also be express (9.14) in the pleasingly symmetrical form

(9.15) $\quad \gamma_1(1-x) - \gamma_1(1+x) = 2\pi \sum_{n=0}^{\infty} \frac{(-1)^{n-1}}{(2n+1)!} \varsigma'(-2n-1)(2\pi x)^{2n+1}$

$$+ [\gamma + \log(2\pi)][\psi(1-x) - \psi(1+x)]$$

With $x = \frac{1}{2}$ we obtain

(9.16) $\quad \sum_{n=0}^{\infty} \frac{(-1)^{n-1}}{(2n+1)!} \varsigma'(-2n-1)(\pi)^{2n+1} = \frac{\gamma + \log \pi}{\pi}$

$\square$

Differentiation of (9.15) results in

(9.17) $\quad -[\gamma_1'(1-x) + \gamma_1'(1+x)] = (2\pi)^2 \sum_{n=0}^{\infty} \frac{(-1)^{n-1}}{(2n)!} \varsigma'(-2n-1)(2\pi x)^{2n}$

$$-[\gamma + \log(2\pi)][\psi'(1+x) + \psi'(1-x)]$$

and with $x = 0$ we have

(9.18) $\quad \gamma_1'(1) = 2\pi^2 \varsigma'(-1) + \varsigma(2)[\gamma + \log(2\pi)]$

It is known that [26]

(9.19) $\quad \gamma_n(x) = \sum_{k=0}^{\infty} \left[ \frac{\log^n(k+x)}{k+x} - \frac{1}{n+1} [\log^{n+1}(k+2) - \log^{n+1}(k+1)] \right]$

Differentiating (9.19) gives us

$$\gamma_n'(x) = \sum_{k=0}^{\infty} \left[ \frac{n \log^{n-1}(k+x) - \log^n(k+x)}{(k+x)^2} \right]$$

which may be expressed as

$$\gamma_n'(x) = (-1)^{n+1} [n \varsigma^{(n-1)}(2,x) + \varsigma^{(n)}(2,x)]$$

In particular we have

(9.20) $\quad \gamma_1'(x) = \varsigma(2,x) + \varsigma'(2,x)$

so that

(9.21) $$\gamma_1'(1) = \varsigma(2) + \varsigma'(2)$$

Using

$$\varsigma'(-1) - \frac{1}{12}(1 - \gamma - \log 2\pi) = \frac{1}{2\pi^2}\varsigma'(2)$$

it is easily seen that (9.18) is the same as (9.21).

Differentiating (7.19) gives us

$$\gamma_1'(x + \tfrac{1}{2}) = 4\gamma_1'(2x) - \gamma_1'(x) + 4\log 2 \cdot \psi'(2x)$$

resulting in

$$\gamma_1'(1) = 4\gamma_1'(1) - \gamma_1'(\tfrac{1}{2}) + 4\log 2 \cdot \psi'(1)$$

and it is easily seen that this concurs with (9.20).

□

There are various series involving the Stieltjes constants and the first derivative $\varsigma'(n)$ at *positive* integer arguments (see for example [20]):

$$\gamma_1(1-x) - \gamma_1 = -\sum_{n=1}^{\infty}\left[\varsigma'(n+1) + H_n^{(1)}\varsigma(n+1)\right]x^n$$

$$\gamma_1(1+x) - \gamma_1 = -\sum_{n=1}^{\infty}(-1)^n\left[\varsigma'(n+1) + H_n^{(1)}\varsigma(n+1)\right]x^n$$

whereupon addition gives us

$$\gamma_1(1-x) + \gamma_1(1+x) - 2\gamma_1 = -2\sum_{n=1}^{\infty}\left[\varsigma'(2n+1) + H_{2n}^{(1)}\varsigma(2n+1)\right]x^{2n}$$

From [24] we have

$$\gamma_0(x) = -\log x - \sum_{n=1}^{\infty}\frac{(-1)^n}{n+1}\varsigma(n+1, x)$$

$$\gamma_1(x) = -\frac{1}{2}\log^2 x + \sum_{n=1}^{\infty}(-1)^n\frac{\varsigma'(n+1, x)}{n+1} + \sum_{n=1}^{\infty}(-1)^n\frac{H_n\varsigma(n+1, x)}{n+1}$$

and these bear some limited structural similarities to (9.15).

**Remark (i)**

Candelpergher, Gadiyar and Padma ([13] and [14, p.68]) have shown that for Re $z > 0$ and $-\pi < \text{Im } z < \pi$

$$(9.22) \qquad \sum_{n=1}^{\infty} e^{-nz} \log n = \sum_{n=0}^{\infty} \frac{(-1)^{n+1}}{n!} \varsigma'(-n) z^n - \frac{\gamma + \log z}{z}$$

and using the relationship

$$\sum_{n=1}^{\infty} (-1)^n a_n = 2 \sum_{n=1}^{\infty} a_{2n} - \sum_{n=1}^{\infty} a_n$$

we see that

$$\sum_{n=1}^{\infty} (-1)^n e^{-nz} \log n = 2 \sum_{n=1}^{\infty} e^{-2nz} \log(2n) - \sum_{n=1}^{\infty} e^{-nz} \log n$$

$$= 2 \log 2 \sum_{n=1}^{\infty} e^{-2nz} + 2 \sum_{n=1}^{\infty} e^{-2nz} \log n - \sum_{n=1}^{\infty} e^{-nz} \log n$$

$$= \frac{2 \log 2}{e^{2z} - 1} + 2 \sum_{n=1}^{\infty} e^{-2nz} \log n - \sum_{n=1}^{\infty} e^{-nz} \log n$$

$$= \frac{2 \log 2}{e^{2z} - 1} + 2 \sum_{n=0}^{\infty} \frac{(-1)^{n+1}}{n!} \varsigma'(-n) 2^n z^n - 2 \frac{\gamma + \log 2z}{2z}$$

$$- \sum_{n=0}^{\infty} \frac{(-1)^{n+1}}{n!} \varsigma'(-n) z^n + \frac{\gamma + \log z}{z}$$

$$= 2 \log 2 \left( \frac{1}{e^{2z} - 1} - \frac{1}{2z} \right) + \sum_{n=0}^{\infty} \frac{(-1)^{n+1}}{n!} \varsigma'(-n) (2^{n+1} - 1) z^n$$

The Bernoulli numbers $B_n$ are given by the generating function

$$\frac{t}{e^t - 1} = \sum_{n=0}^{\infty} B_n \frac{t^n}{n!} \qquad , (|t| < 2\pi)$$

Therefore, employing $B_0 = 1$, $B_1 = -\frac{1}{2}$ and $B_{2n+1} = 0$ for all $n \geq 1$, we have

$$\frac{1}{e^{2z} - 1} = \sum_{n=0}^{\infty} 2^{n-1} B_n \frac{z^{n-1}}{n!}$$

$$= \frac{1}{2z} - \frac{1}{2} + \sum_{n=2}^{\infty} 2^{n-1} B_n \frac{z^{n-1}}{n!}$$

$$= \frac{1}{2z} - \frac{1}{2} + \sum_{n=1}^{\infty} 2^{2n-1} B_{2n} \frac{z^{2n-1}}{(2n)!}$$

Hence, we deduce that

(9.23) $\quad \sum_{n=1}^{\infty} (-1)^n e^{-nz} \log n = 2\log 2 \left( -\frac{1}{2} + \sum_{n=1}^{\infty} 2^{2n-1} B_{2n} \frac{z^{2n-1}}{(2n)!} \right) + \sum_{n=0}^{\infty} \frac{(-1)^{n+1}}{n!} \varsigma'(-n)(2^{n+1}-1) z^n$

We note that the right-hand side is finite at $z = 0$ and, in the limit as $z \to 0$, we see that

(9.24) $\quad \lim_{z \to 0} \left[ \sum_{n=1}^{\infty} (-1)^n e^{-nz} \log n \right] = \varsigma'(0) - \log 2$

Differentiating (9.23) results in

(9.25) $\quad -\sum_{n=1}^{\infty} (-1)^n n e^{-nz} \log n = 2\log 2 \sum_{n=1}^{\infty} 2^{2n-1} (2n-1) B_{2n} \frac{z^{2n-2}}{(2n)!} + \sum_{n=0}^{\infty} \frac{(-1)^{n+1}}{n!} \varsigma'(-n)(2^{n+1}-1) n z^{n-1}$

and we have

(9.26) $\quad \lim_{z \to 0} \left[ \sum_{n=1}^{\infty} (-1)^n n e^{-nz} \log n \right] = -\frac{1}{3} \log 2 - 3\varsigma'(-1)$

The second derivative of (9.23) gives us

$$\sum_{n=1}^{\infty} (-1)^n n^2 e^{-nz} \log n = 2\log 2 \sum_{n=1}^{\infty} 2^{2n-1} B_{2n} (2n-1)(2n-2) \frac{z^{2n-3}}{(2n)!}$$

$$+ \sum_{n=0}^{\infty} \frac{(-1)^{n+1}}{n!} \varsigma'(-n)(2^{n+1}-1) n(n-1) z^{n-2}$$

In the limit as $z \to 0$ we see that

$$-7\varsigma'(-2) = \lim_{z \to 0} \left[ \sum_{n=1}^{\infty} (-1)^n n^2 e^{-nz} \log n \right]$$

We have

$$\varsigma'(-2) = -\frac{\varsigma(3)}{4\pi^2}$$

The alternating Riemann zeta function is defined by $\varsigma_a(s) = \sum_{n=1}^{\infty} \frac{(-1)^{n+1}}{n^s}$ and we have

$$\varsigma_a(s) = (1 - 2^{1-s}) \varsigma(s)$$

Hence we see that $\varsigma_a(3) = \frac{3}{4}\varsigma(3)$ and we therefore obtain

$$(9.27) \qquad \varsigma_a(3) = \frac{3}{7}\pi^2 \lim_{z \to 0}\left[\sum_{n=1}^{\infty}(-1)^n n^2 e^{-nz} \log n\right]$$

which corresponds with the Euler summation given by Candelpergher [14, p.146].

This result is very reminiscent of Euler's work as reported by Ayoub [6]. Euler stated that

$$\varsigma_a(3) = \frac{3}{7}\pi^2 \sum_{n=1}^{\infty}(-1)^n n^2 \log n$$

However (and it's a *big* however), it should be noted that Euler's series is obviously not convergent.

It may be observed that

$$\sum_{n=1}^{\infty}(-1)^n n^2 e^{-nz} \log n = \frac{\partial}{\partial s}\sum_{n=1}^{\infty}(-1)^{n+1}\frac{e^{-nz}}{n^s}\bigg|_{s=-2}$$

$$= -\frac{\partial}{\partial s} Li_s(-e^{-z})\bigg|_{s=-2}$$

in terms of the polylogarithm function.

$\square$

Integrating (9.23) gives us

$$(9.28) \qquad \sum_{n=1}^{\infty}(-1)^n e^{-nz}\frac{\log n}{n} - \sum_{n=1}^{\infty}(-1)^n\frac{\log n}{n} = 2\log 2\left(-\frac{1}{2}z + \sum_{n=1}^{\infty}2^{2n-1}B_{2n}\frac{z^{2n}}{2n(2n)!}\right)$$

$$+ \sum_{n=0}^{\infty}\frac{(-1)^{n+1}}{(n+1).n!}\varsigma'(-n)(2^{n+1}-1)z^{n+1}$$

**Remark (ii)**

Letting $z = \alpha + 2\pi i x$ in (9.22) with $-\pi < \operatorname{Im} z < \pi$ gives us for the imaginary part

$$(9.29) \qquad \lim_{\alpha \to 0}\sum_{n=1}^{\infty} e^{-\alpha n}\log n \cdot \sin 2\pi nx = \sum_{n=0}^{\infty}\frac{(-1)^{n-1}}{(2n+1)!}\varsigma'(-2n-1)(2\pi x)^{2n+1} - \frac{\gamma + \log(2\pi x)}{2\pi x}$$

and we note that this has the same form as (9.12) below

$$\lim_{s\to 0}\sum_{n\geq 1}\log n\,\frac{\sin 2\pi nx}{n^s}=\sum_{n=0}^{\infty}\frac{(-1)^{n-1}}{(2n+1)!}\varsigma'(-2n-1)(2\pi x)^{2n+1}-\frac{\gamma+\log(2\pi x)}{2\pi x}$$

The real part gives us

(9.30) $$\lim_{\alpha\to 0}\sum_{n=1}^{\infty}e^{-\alpha n}\log n\cdot\cos 2\pi nx=-\sum_{n=0}^{\infty}\frac{(-1)^n}{(2n)!}\varsigma'(-2n)(2\pi x)^{2n}-\frac{1}{4x}$$

which corresponds with (9.9).

**Remark (iii)**

A *very loose* connection between $\sum_{n=1}^{\infty}e^{-nz}\log n$ and $\sum_{n=0}^{\infty}\frac{(-1)^{n+1}}{n!}\varsigma'(-n)z^n$ may be formally seen below:

$$\sum_{n=1}^{\infty}e^{-nz}\log n=\sum_{n=1}^{\infty}\sum_{k=0}^{\infty}\frac{(-1)^k z^k}{k!}n^k\log n$$

$$\sim\sum_{k=0}^{\infty}\frac{(-1)^k z^k}{k!}\sum_{n=1}^{\infty}n^k\log n$$

$$\sim\sum_{k=0}^{\infty}\frac{(-1)^{k-1}\varsigma'(-k)}{k!}z^k$$

## 10. Open access to our own work

This paper contains references to a number of other papers and, fortunately, most of them are currently freely available on the internet. Surely now is the time that all of <u>our</u> work should be freely accessible by <u>all</u>. The mathematics community should lead the way on this by publishing <u>everything</u> on arXiv, or in an equivalent open access repository. We think it, we write it, so why hide it? You know it makes sense.

### ACKNOWLEDGEMENTS

I am greatly indebted to Bernard Candelpergher for kindly providing me with the substantive part of Section 9.

Wessex House,
Devizes Road,
Upavon,
Pewsey,
Wiltshire SN9 6DL